\documentclass[11pt]{article}

\usepackage{hyperref}
\usepackage{rate-distortion}
\usepackage{fullpage}
\usepackage[numbers]{natbib}
\usepackage{color}
\usepackage{bbm}
\usepackage{psfrag}

\definecolor{darkblue}{rgb}{.15,0,.7}
\hypersetup{
  colorlinks=true,
  linkcolor=darkblue,
  citecolor=darkblue,
  urlcolor=darkblue
}
\renewcommand{\onevec}{\mathbbm{1}}

\newcommand{\minimax}{\mathfrak{M}}
\newcommand{\parameterspace}{\Theta}
\newcommand{\laplace}{\mathop{\rm Laplace}}

\newcommand{\project}{\Pi}
\newcommand{\metric}{\rho}
\newcommand{\numobs}{n}
\newcommand{\QFAMA}{\mc{\channelprob}_\diffp} 

\newcommand{\packset}{\ensuremath{\mc{V}}}  
\newcommand{\packrv}{V}   
\newcommand{\packval}{\nu}
\newcommand{\altpackval}{\nu'}
\newcommand{\meanstatprob}{\overline{\statprob}}

\newcommand{\channeldomain}{\mc{\channelrv}}
\newcommand{\linfset}{\mc{G}_\diffp} 
\newcommand{\linfsetfinite}[1]{\mc{G}_{\diffp, {#1}}} 
\newcommand{\optdens}{\gamma} 

\newcommand{\characteristic}[1]{\mathsf{1}_{#1}} 
\newcommand{\histelement}{\mathsf{e}} 

\newcommand{\simplex}{\Delta} 

\newcommand{\numderiv}{\beta} 
\newcommand{\lebesgue}{\mu} 
\newcommand{\basisfunc}{\varphi} 
\newcommand{\lipconst}{C} 
\newcommand{\numbin}{\ensuremath{k}} 
\newcommand{\sbound}{\ensuremath{B}} 
\newcommand{\orthbound}{\ensuremath{B_0}} 
\newcommand{\bernoulli}{\mathop{\rm Bernoulli}} 

\newcommand{\test}{\psi}

\definecolor{innerboxcolor}{rgb}{.9,.95,1}
\definecolor{outerlinecolor}{rgb}{.6,0,.2}

\newcommand{\densclass}[1][]{
  \ifthenelse{\isempty{#1}}{%
    \mc{F}_\numderiv
  }{
    \mc{F}_{\numderiv, {#1}}
  }
}

\makeatletter
\long\def\@makecaption#1#2{
        \vskip 0.8ex
        \setbox\@tempboxa\hbox{\small {\bf #1:} #2}
        \parindent 1.5em  
        \dimen0=\hsize
        \advance\dimen0 by -3em
        \ifdim \wd\@tempboxa >\dimen0
                \hbox to \hsize{
                        \parindent 0em
                        \hfil 
                        \parbox{\dimen0}{\def\baselinestretch{0.96}\small
                                {\bf #1.} #2
                                } 
                        \hfil}
        \else \hbox to \hsize{\hfil \box\@tempboxa \hfil}
        \fi
        }
\makeatother

\long\def\comment#1{}

\newcommand{\nnz}{\ensuremath{k}} 

\newcommand{\locnnz}{\ensuremath{\ell}} 

\newcommand{\Mpack}{\ensuremath{K}}

\newcommand{\USET}{\ensuremath{\mathcal{U}_\Mpack}}

\newcommand{\MySet}{\ensuremath{D}}

\begin{document}

\begin{center}
  {\Large \textbf{Local Privacy and Minimax Bounds: \\
      Sharp Rates for Probability Estimation}}

  \vspace{1cm}
  
  \begin{tabular}{ccccc}
    John C.\ Duchi$^\dagger$ & & Michael I.\ Jordan$^{
      \dagger, \ast}$ &  & Martin
    J.\ Wainwright$^{\dagger, \ast}$
    \\ \texttt{jduchi@eecs.berkeley.edu} & &
    \texttt{jordan@stat.berkeley.edu} & &
    \texttt{wainwrig@stat.berkeley.edu} \\
  \end{tabular}

  \vspace{.5cm}
  \begin{tabular}{ccc}
    Department of Statistics$^\ast$ & & Department of Electrical
    Engineering and Computer Science$^{\dagger}$
  \end{tabular}

  \vspace{.2cm}
  University of California, Berkeley \\
  Berkeley, CA, 94720
\end{center}

\begin{abstract}
  We provide a detailed study of the estimation of probability
  distributions---discrete and continuous---in a stringent setting 
  in which data is kept private even from the statistician.  We
  give sharp minimax rates of convergence for estimation in these
  locally private settings, exhibiting fundamental tradeoffs between
  privacy and convergence rate, as well as providing tools to allow
  movement along the privacy-statistical efficiency continuum. One of the
  consequences of our results is that Warner's classical work on
  randomized response is an optimal way to perform survey sampling
  while maintaining privacy of the respondents.
\end{abstract}


\section{Introduction}

The original motivation for providing privacy in statistical problems,
first discussed by~\citet{Warner65}, was that ``for reasons of
modesty, fear of being thought bigoted, or merely a reluctance to
confide secrets to strangers,'' respondents to surveys might prefer to
be able to answer certain questions non-truthfully, or at least
without the interviewer knowing their true response.  With this
motivation, Warner considered the problem of estimating the fractions
of the population belonging to certain strata, which can be viewed as
probability estimation within a multinomial model.  In this paper, we
revisit Warner's probability estimation problem, doing so within a 
theoretical framework that allows us to characterize optimal
estimation under constraints on privacy.  We also apply our theoretical
tools to a further probability estimation problem---that of nonparametric 
density estimation.

In the large body of research on privacy and statistical
inference~\cite[e.g.,][]{Warner65, 
Fellegi72, DuncanLa86, DuncanLa89, FienbergMaSt98}, a major focus has been
on the problem of reducing disclosure risk: the probability that a member 
of a dataset can be identified given released statistics of the dataset.
The literature has stopped short, however, of providing a formal treatment 
of disclosure risk that would permit decision-theoretic tools to be used
in characterizing tradeoffs between the utility of achieving privacy and the
utility associated with an inferential goal.  Recently, a formal treatment 
of disclosure risk known as ``differential privacy'' has been proposed 
and studied in the cryptography, database and theoretical computer science
literatures~\cite{DworkMcNiSm06, BarakChDwKaMcTa07, Dwork08}. Differential 
privacy
has strong semantic privacy guarantees make it a good candidate for declaring 
a statistical procedure or data collection mechanism private, and
has been the focus of a growing body of recent work~\cite{
  Dwork08, EvfimievskiGeSr03, HardtTa10, WassermanZh10,
  Smith11, ChaudhuriMoSa11,
  KasiviswanathanLeNiRaSm11, De12, ChaudhuriHs12, DuchiJoWa13_parametric}.

In this paper, we bring together the formal treatment of disclosure risk 
provided by differential privacy with the tools of minimax decision theory
to provide a theoretical treatment of probability estimation under privacy
constraints.  Just as in classical minimax theory, we are able to provide 
lower bounds on the convergence rates of any estimator, in our case under 
a restriction to estimators that guarantee privacy.  We complement these 
results with matching upper bounds that are achievable using computationally 
efficient algorithms.  We thus bring classical notions of privacy, as introduced 
by \citet{Warner65}, into contact with differential privacy and statistical 
decision theory, obtaining quantitative tradeoffs between privacy and statistical 
efficiency.

\subsection{Setting and contributions}

Let us develop a bit of basic formalism before describing---at a high
level---our main results.  We study procedures that receive private views
$\channelrv_1, \ldots, \channelrv_n \in \channeldomain$ of an original set of
observations, $\statrv_1, \ldots, \statrv_n \in \statdomain$, where $\statdomain$ is
the (known) sample space.  In our setting, $\channelrv_i$ is drawn conditional
on $\statrv_i$ via the \emph{channel distribution} $\channelprob(\channelrv_i
\mid \statrv_i = \statsample, \channelrv_j = \channelval_j, j \neq i)$.  Note
that this channel allows ``interactivity''~\cite{DworkMcNiSm06}, meaning that
the distribution of $\channelrv_i$ may depend on $\statrv_i$ as well as the
private views $\channelrv_j$ of $\statrv_j$ for $j \neq i$. Allowing
interactivity---rather than forcing $\channelrv_i$ to be independent of
$\channelrv_j$---in some cases allows more efficient algorithms, and in our
setting means that our lower bounds are stronger.

We assume each of these private views $\channelrv_i$ is
$\diffp$-differentially private for the original data $\statrv_i$. To give a
precise definition for this type of privacy, known as ``local privacy,'' 
let $\sigma(\channeldomain)$ be the $\sigma$-field on $\channeldomain$ over 
which the channel $\channelprob$ is defined.  Then $\channelprob$ provides 
\emph{$\diffp$-local differential privacy} if
\begin{equation}
  \sup \left\{
  \frac{\channelprob(S \mid \statrv_i = \statsample,
    \channelrv_j = \channelval_j, j \neq i)}{
    \channelprob(S \mid \statrv_i = \statsample',
    \channelrv_j = \channelval_j, j \neq i)}
  \mid
  S \in \sigma(\channeldomain),
  \channelval_j \in \channeldomain, ~ \mbox{and} ~
  \statsample, \statsample' \in \statdomain \right\}
  \le \exp(\diffp).
  \label{eqn:local-privacy}
\end{equation}
In the non-interactive setting (in which we impose the constraint that 
the providers of the data release a private view independently of the 
other data providers) the expression~\eqref{eqn:local-privacy} simplifies to
\begin{equation}
  \sup_{S \in \sigma(\channeldomain)}
  \sup_{\statsample, \statsample' \in \statdomain}
  \frac{\channelprob(S \mid \statrv = \statsample)}{
    \channelprob(S \mid \statrv = \statsample')}
  \le \exp(\diffp),
  \label{eqn:local-privacy-no-interaction}
\end{equation}
a formulation of local privacy first proposed by 
\citet{EvfimievskiGeSr03}. Although more complex to analyze, the likelihood ratio
bound~\eqref{eqn:local-privacy} is attractive for many reasons. It means that
any individual providing data guarantees his or her own privacy---no further
processing or mistakes by a collection agency can compromise one's data---and 
the individual has plausible deniability about taking a value $\statsample$,
since any outcome $\channelval$ is nearly as likely to have come from some 
other initial value $\statsample'$. The likelihood ratio also controls 
the error rate in tests for the presence of points $\statsample$ in
the data~\cite{WassermanZh10}. All that is required is that the likelihood
ratio~\eqref{eqn:local-privacy} be bounded no matter the data provided by
other participants.

In the current paper, we study minimax convergence rates when the data 
provided satisfies the local privacy guarantee~\eqref{eqn:local-privacy}.  
Our two main results quantify the penalty that must be paid when local 
privacy at a level $\diffp$ is provided in multinomial estimation and 
density estimation problems.  At a high level, our first result implies 
that for estimation of a $d$-dimensional multinomial probability mass 
function, the effective sample size of \emph{any} statistical estimation 
procedure decreases from $n$ to $n \diffp^2 / d$ whenever $\diffp$
is a sufficiently small constant. A consequence of our results
is that Warner's randomized response procedure~\cite{Warner65}
enjoys optimal sample complexity; it is interesting to note that
even with the recent focus on privacy and statistical inference, 
the optimal privacy-preserving strategy for problems such as survey 
collection has been known for almost 50 years.

Our second main result, on density estimation, exhibits an interesting 
departure from standard minimax estimation results. If the density 
being estimated has $\numderiv$ continuous derivatives, then classical 
results on density estimation \cite[e.g.,][]{Yu97, YangBa99, Tsybakov09}) 
show that the minimax integrated squared error scales (in the number 
$n$ of samples) as $n^{-2 \numderiv / (2 \numderiv + 1)}$. 
In the locally private case, we show that---even when
densities are bounded and well-behaved---there is a difference in the
\emph{polynomial} rate of convergence: we obtain a scaling of
$(\diffp^2 n)^{-2 \numderiv / (2 \numderiv + 2)}$.  We give
efficiently implementable algorithms
that attain sharp upper bounds as companions to our lower
bounds, which in some cases exhibit the necessity of non-trivial
sampling strategies to guarantee privacy.

\paragraph{Notation:}
We summarize here the notation used
throughout the paper. Given distributions $P$ and $Q$ defined on a
space $\statdomain$, each absolutely continuous with respect to a
distribution $\mu$ (with corresponding densities $p$ and $q$), the
KL-divergence between $P$ and $Q$ is defined by
\begin{equation*}
  \dkl{P}{Q} \defeq \int_\statdomain dP \log \frac{dP}{dQ} =
  \int_\statdomain p \log \frac{p}{q} d\mu.
\end{equation*}
Letting $\sigma(\statdomain)$ denote the (an appropriate)
$\sigma$-field on $\statdomain$, the total variation distance between
the distributions $P$ and $Q$ on $\statdomain$ is given by
\begin{equation*}
  \tvnorm{P - Q} \defeq \sup_{S \in \sigma(\statdomain)} |P(S) - Q(S)|
  = \half \int_\statdomain \left|p(\statsample) -
  q(\statsample)\right| d\mu(\statsample).
\end{equation*}
For random vectors $X$ and $Y$, where $X$ is distributed according to the
distribution $P$ and $Y \mid X$ is distributed according to $Q(\cdot \mid X)$,
let $\marginprob(\cdot) = \int Q(\cdot \mid x) dP(x)$ denote the marginal
distribution of $Y$. The mutual information between $X$ and $Y$ is
\begin{align*}
\information(X; Y) \defeq \E_P \left[ \dkl{ Q(\cdot \mid
    X)}{\marginprob(\cdot)} \right] = \int \dkl{ Q(\cdot \mid X = x)}
            {\marginprob(\cdot)} d P(x).
\end{align*}
A random variable $Y$ has $\laplace(\alpha)$ distribution if its
density $p_Y(y) = \frac{\alpha}{2} \exp\left(-\alpha
|y|\right)$, where $\alpha > 0$.  For matrices $A, B \in \R^{d \times
  d}$, the notation $A \preceq B$ means $B - A$ is positive
semidefinite, and $A \prec B$ means $B - A$ is positive definite. We
write $a_n \lesssim b_n$ to denote that $a_n = \order(b_n)$ and $a_n
\asymp b_n$ to denote that $a_n = \order(b_n)$ and $b_n =
\order(a_n)$.  For a convex set $C \subset \R^d$, we let $\project_C$
denote the orthogonal projection operator onto $C$,
i.e., $\project_C(v) \defeq \argmin_{w \in C} \{\ltwo{v - w}\}$.



\section{Background and Problem Formulation}

In this section, we provide the necessary background on the minimax
framework used throughout the paper.  Further details on minimax
techniques can be found in several standard sources~\cite[e.g.,][]{Birge83, 
YangBa99, Yu97, Tsybakov09}.  We also reference our companion 
paper on parametric statistical inference under differential privacy
constraints~\cite{DuchiJoWa13_parametric}; we make use of
two theorems from that earlier paper, but in order to keep the current 
paper self-contained, we restate them in this section.

\subsection{Minimax framework}
\label{sec:minimax-framework}

Let $\mc{P}$ denote a class of distributions on the sample space
$\statdomain$, and let $\optvar : \mc{P} \rightarrow \parameterspace$
denote a function defined on $\mc{P}$. The range 
$\parameterspace$ depends on the underlying statistical
model; for example, for density estimation, $\parameterspace$ may
consist of the set of probability densities defined on $[0, 1]$.  We
let $\metric$ denote the semi-metric on the space $\parameterspace$
that we use to measure the error of an estimator for $\optvar$, and
$\Phi : \R_+ \rightarrow \R_+$ be a non-decreasing function with
$\Phi(0) = 0$ (for example, $\Phi(t) = t^2$).

Recalling that $\channeldomain$ is the domain of the private variables
$\channelrv_i$, let $\what{\optvar} : \channeldomain^n \rightarrow
\optdomain$ denote an arbitrary estimator for $\optvar$.  Let
$\QFAMA$ denote the set of conditional (or channel) distributions
guaranteeing $\diffp$-local privacy~\eqref{eqn:local-privacy}; then for any
$Q \in \QFAMA$ we can define the minimax rate
\begin{subequations}
  \begin{equation}
    \minimax_n\left(\optvar(\mc{P}), \Phi \circ \metric, \channelprob\right)
    \defeq \inf_{\what{\optvar}} \sup_{\statprob \in \mc{P}}
    \E_{\statprob, \channelprob}
    \left[\Phi\left(\metric(\what{\optvar}(\channelrv_1,
      \ldots, \channelrv_n), \optvar(\statprob))\right)\right]
    \label{eqn:minimax-risk-Q}
  \end{equation}
associated with estimating $\optvar$ based on the private samples
$(\channelrv_1, \ldots, \channelrv_\numobs)$.  In the
definition~\eqref{eqn:minimax-risk-Q}, the expectation is taken both
with respect to the distribution $\statprob$ on the variables
$\statrv_1, \ldots, \statrv_n$ and the $\diffp$-private channel
$\channelprob$.  By taking the infimum over all possible channels $Q
\in \QFAMA$, we obtain the central object of interest for this
paper, the \emph{$\diffp$-private minimax rate} for the family
$\theta(\mc{P})$, defined as
  \begin{equation}
    \label{eqn:minimax-risk}
    \minimax_n(\optvar(\mc{P}), \Phi \circ \metric, \diffp) \defeq
    \inf_{\what{\theta}, \channelprob \in \QFAMA} \sup_{\statprob \in
      \mc{P}} \E_{\statprob, \channelprob}
    \left[\Phi\left(\metric(\what{\optvar}(\channelrv_1, \ldots,
      \channelrv_n), \optvar(\statprob))\right)\right].
  \end{equation}
\end{subequations}

A standard route for lower bounding the minimax
risk~\eqref{eqn:minimax-risk-Q} is by reducing the estimation problem
to the testing problem of identifying a point $\optvar \in \optdomain$
from a finite collection of well-separated points~\cite{Yu97,
  YangBa99}. Given an index set $\packset$ of finite cardinality, the
indexed family of distributions $\{\statprob_\packval, \packval \in
\packset\} \subset \mc{P}$ is said to be a $2\delta$-packing of
$\optdomain$ if $\metric(\optvar(\statprob_\packval),
\optvar(\statprob_{\altpackval})) \ge 2 \delta$ for all $\packval \neq
\altpackval$ in $\packset$. The setup is that of a standard hypothesis testing
problem: nature chooses $\packrv \in \packset$ uniformly at random,
then data $(\statrv_1, \ldots, \statrv_\numobs)$ are drawn from the 
$\numobs$-fold conditional product distribution $\statprob_\packval^n$, 
conditioning on $\packrv = \packval$.  The problem is to identify the 
member $\packval$ of the packing set $\packset$.

In this work we have the additional complication that all the statistician
observes are the private samples
$\channelrv_1, \ldots, \channelrv_\numobs$.
To that end, if we let
$\channelprob^n(\cdot \mid \statsample_{1:n})$ denote the conditional
distribution of $\channelrv_1, \ldots, \channelrv_n$ given that $\statrv_1 =
\statsample_1, \ldots, \statrv_n = \statsample_n$, we define the marginal
channel $\marginprob^n_\packval$ via the expression
\begin{equation}
  \marginprob_\packval^n(A) \defeq
  \int \channelprob^n(A \mid \statsample_1, \ldots, \statsample_n)
  d \statprob_\packval(\statsample_1, \ldots, \statsample_n)
  ~~~ \mbox{for}~ A \in \sigma(\channeldomain^n).
  \label{eqn:marginal-channel}
\end{equation}
Letting $\test : \channeldomain^n \rightarrow \packset$ denote an arbitrary
testing procedure---a measureable mapping $\channeldomain^n \to \packset$---we
have the following minimax risk bound, whose two parts are known as
Le Cam's two-point method and Fano's inequality. In the lemma, we let $\P$
denote the joint distribution of the random variable $\packrv$ and the samples
$\channelrv_i$.
\begin{lemma}[Minimax risk bound]
  \label{lemma:minimax-risk-bound}
  For the previously described estimation and testing problems, we have
  the lower bound
  \begin{align}
    \minimax_n(\optvar(\mc{P}), \Phi \circ \metric, \channelprob)
    \ge \Phi(\delta) \, \inf_{\test}
    \P(\test(\channelrv_1, \ldots, \channelrv_n) \neq \packrv),
  \end{align}
  where the infimum is taken over all testing procedures. For
  a binary test specified by $\packset = \{\packval, \altpackval\}$,
  \begin{subequations}
    \begin{align}
      \label{eqn:le-cam}
      \inf_{\test} \P \left( \test(\channelrv_1, \ldots, \channelrv_n) \neq
      \packrv \right ) = \half - \half \tvnorm{\marginprob_\packval^n -
        \marginprob_{\altpackval}^n},
    \end{align}
    and more generally,
    \begin{align}
      \label{eqn:fano}
      \inf_{\test} \P(\test(\channelrv_1, \ldots, \channelrv_n) \neq
      \packrv) \ge \left [1 - \frac{\information(\channelrv_1, \ldots,
          \channelrv_n; \packrv) + \log 2}{\log |\packset|} \right].
    \end{align}
  \end{subequations}
\end{lemma}
\noindent
For Le Cam's inequality~\eqref{eqn:le-cam}, see, e.g., Lemma~1
of~\citet{Yu97} or Theorem 2.2 of~\citet{Tsybakov09}; for Fano's
inequality~\eqref{eqn:fano}, see Eq.~(1)
of~\citet{YangBa99} or Chapter~2 of~\citet{CoverTh06}.

\subsection{Information bounds}

The main step in proving minimax lower bounds is to control
the divergences involved in the lower bounds~\eqref{eqn:le-cam}
and~\eqref{eqn:fano}. In our companion paper~\cite{DuchiJoWa13_parametric}, 
we present two results, which we now review, in which bounds on 
$\tvnorm{\marginprob_\packval^n - \marginprob_{\altpackval}^n}$ 
and $\information(\channelrv_1, \ldots, \channelrv_n; \packrv)$ 
are obtained as a function of the amount of privacy provided 
and the distances between the underlying distributions
$\statprob_\packval$.  The first result~\cite[Theorem 1 and 
Corollary 1]{DuchiJoWa13_parametric} gives control over pairwise 
KL-divergences between the marginals~\eqref{eqn:marginal-channel}, 
which lends itself to application of Le Cam's method~\eqref{eqn:le-cam} 
and simple applications of Fano's inequality~\eqref{eqn:fano}. The second
result~\cite[Theorem 2 and Corollary 4]{DuchiJoWa13_parametric}
provides a variational upper bound on the mutual information
$\information(\channelrv_1, \ldots, \channelrv_\numobs;
\packrv)$---variational in the sense that it requires optimization
over the set of functions
\begin{equation*}
  \linfset \defeq \Big\{\optdens \in L^\infty(\statdomain) \, \mid \,
  \sup_{\statsample \in \statdomain} |\optdens(\statsample)| \le \half
  (e^\diffp - e^{-\diffp}) \Big\}.
\end{equation*}
Here \mbox{$L^\infty(\statdomain) \defeq \{f : \statdomain \rightarrow
  \R \mid \sup_{\statsample \in \statdomain} |f(\statsample)| <
  \infty\}$} denotes the space of bounded functions on $\statdomain$.

Our bounds apply to any channel distribution $\channelprob$ that is
$\diffp$-locally private~\eqref{eqn:local-privacy}.  For each \mbox{$i
  \in \{1, \ldots, \numobs\}$,} let $\statprob_{\packval, i}$ be the
distribution of $\statrv_i$ conditional on the random packing element
$\packrv = \packval$, and let $\marginprob_\packval^n$ be the marginal
distribution~\eqref{eqn:marginal-channel} induced by passing
$\statrv_i$ through $\channelprob$.  Define the mixture distribution
\mbox{$\meanstatprob_i = \frac{1}{|\packset|} \sum_{\packval \in
    \packset} \statprob_{\packval, i}$,} and the linear functionals
$\varphi_{\packval, i} : L^\infty(\statdomain) \rightarrow \R$ by
  \begin{equation*}
    \varphi_{\packval, i}(\optdens)
    \defeq \int_{\statdomain} \optdens(\statsample) \left(
    d\statprob_{\packval, i}(\statsample) - d\meanstatprob_i(\statsample)
    \right).
  \end{equation*}
With this notation we can state the following proposition, which summarizes 
the results that we will need from~\citet{DuchiJoWa13_parametric}:
\begin{proposition}[Information bounds]
  \label{proposition:information-bounds}
\begin{itemize}
\item[(a)] For all $\diffp \ge 0$,
  \begin{equation}
    \label{eqn:one-d-kl}
    \dkl{\marginprob_\packval^n}{\marginprob_{\altpackval}^n} +
    \dkl{\marginprob_{\altpackval}^n}{\marginprob_\packval^n} \le 4
    (e^\diffp - 1)^2 \sum_{i = 1}^n \tvnorm{\statprob_{\packval, i} -
      \statprob_{\altpackval, i}}^2.
  \end{equation}
\item[(b)] For all $\diffp \in [0, \log(\half + \half \sqrt{3})]$,
  \begin{equation}
    \label{eqn:super-fano}
    \information(\channelrv_1, \ldots, \channelrv_n; \packrv)
    \le C_\diffp \sum_{i=1}^n \frac{1}{|\packset|}
    \sup_{\optdens \in \linfset} \sum_{\packval \in \packset}
    \left(\varphi_{\packval,i}(\optdens)\right)^2,
  \end{equation}
where $C_\diffp \defn 4 / (e^{-\diffp} - 2 (e^\diffp - 1))$.
\end{itemize}
\end{proposition}

By combining Proposition~\ref{proposition:information-bounds} with
Lemma~\ref{lemma:minimax-risk-bound}, it is possible to derive sharp
lower bounds on arbitrary estimation procedures under $\diffp$-local
privacy.  In particular, we may apply the bound~\eqref{eqn:one-d-kl}
with Le Cam's method~\eqref{eqn:le-cam}, though lower bounds so
obtained often lack dimension dependence we might hope to capture (see
Section~3.2 of~\citet{DuchiJoWa13_parametric} for more discussion of
this issue).  The bound~\eqref{eqn:super-fano}, which (up to constants) 
implies the bound~\eqref{eqn:one-d-kl}, allows more careful control 
via suitably constructed packing sets $\packset$ and application 
of Fano's inequality~\eqref{eqn:fano}, since the supremum
controls a more global view of the structure of $\packset$.  In the
rest of this paper, we demonstrate this combination for probability 
estimation problems.


\section{Multinomial Estimation under Local Privacy}
\label{sec:multinomial-estimation}

In this section we return to the classical problem of avoiding answer bias in
surveys, the original motivation for studying local privacy~\cite{Warner65}.
We provide a detailed study of estimation of a multinomial probability under
$\diffp$-local differential privacy, providing sharp upper and lower bounds
for the minimax rate of convergence.

\subsection{Minimax rates of convergence for multinomial estimation}

Consider the
probability simplex \mbox{$\simplex_d \defeq \big \{\optvar \in \R^d
  \, \mid \, \optvar \ge 0, \sum_{j=1}^d \optvar_j = 1 \big \}$} in
$\R^d$. The multinomial estimation problem is defined as follows.
Given a vector $\optvar \in \simplex_d$, samples $\statrv$
are drawn i.i.d.\ from a multinomial with parameters $\optvar$, where
$\statprob_\optvar(\statrv = j) = \optvar_j$ for $j \in \{1, \ldots,
d\}$, and our goal is to estimate the probability vector $\optvar$.
In one of the earliest evaluations of privacy, \citet{Warner65}
studied the Bernoulli variant of this problem, and proposed a simple
privacy-preserving mechanism known as \emph{randomized response}: 
for a given survey question, respondents provide a truthful answer 
with probability $p > 1/2$, and a lie with probability $1 -p$.

In our setting, we assume that the statistician sees random variables
$\channelrv_i$ that are all $\diffp$-locally
private~\eqref{eqn:local-privacy} for the corresponding samples
$\statrv_i$ from the multinomial. In this case, we have the following
result, which characterizes the minimax rate of estimation of a
multinomial in terms of mean-squared error $\E[\ltwos{\what{\optvar} -
    \optvar}^2]$.
\begin{theorem}
  \label{theorem:multinomial-rate}
There exist universal constants $0 < c_\ell \le c_u < 5$ such that for
all $\diffp \in [0, 1/4]$, the minimax rate for multinomial estimation
satisfies the bounds
\begin{equation}
\label{eqn:multinomial-rate}
 c_\ell \max_{k \in \{1, \ldots, d\}} \min \left\{ \frac{1}{k},
 \frac{k \log \frac{d}{k}}{n \diffp^2} \right\} \; \leq \; \minimax_n
 \left( \simplex_d, \ltwo{\cdot}^2, \diffp\right) \; \le \; c_u \min
 \left \{1, \frac{d}{n \diffp^2} \right\}.
  \end{equation}
\end{theorem}
\noindent
We provide a proof of the lower bound in
Theorem~\ref{theorem:multinomial-rate} in
Section~\ref{sec:proof-multinomial-rate}. Simple estimation strategies achieve
the lower bound, and we believe exploring them is interesting, so we provide
them in the next section.

Theorem~\ref{theorem:multinomial-rate} shows that providing local
privacy can sometimes be quite detrimental to the quality of statistical
estimators. Indeed, let us compare this rate to the classical rate in
which there is no privacy. Then estimating $\optvar$ via proportions 
(i.e., maximum likelihood), we have
\begin{equation*}
  \E\left[\ltwos{\what{\optvar} - \optvar}^2\right]
  = \sum_{j=1}^d \E\left[(\what{\optvar}_j - \optvar_j)^2\right]
  = \frac{1}{n} \sum_{j=1}^d \optvar_j(1 - \optvar_j)
  \le \frac{1}{n}\left(1 - \frac{1}{d}\right)
  < \frac{1}{n}.
\end{equation*}
On the other hand, an appropriate choice of $\nnz$ in
Theorem~\ref{theorem:multinomial-rate} implies that
\begin{equation}
  \label{eqn:multinomial-limits}
  \min\left\{1, \frac{1}{\sqrt{n \diffp^2}},
  \frac{d}{n \diffp^2}\right\}
  \lesssim \minimax_n\left(\simplex_d, \ltwo{\cdot}^2, \diffp\right)
  \lesssim \min\left\{1, \frac{d}{n \diffp^2}\right\},
\end{equation}
which we show in Section~\ref{sec:proof-multinomial-rate}.
Thus, for suitably large sample sizes $n$, the effect of providing
differential privacy at a level $\diffp$ causes a reduction in the effective
sample size of $n \mapsto n \diffp^2 / d$.

\subsection{Private multinomial estimation strategies}
\label{sec:private-multinomial-estimation}

An interesting consequence of the lower bound in~\eqref{eqn:multinomial-rate} 
is the following fact that we now demonstrate: Warner's classical 
randomized response mechanism~\cite{Warner65} (with minor
modification) achieves the optimal convergence rate. Thus, although it
was not originally recognized as such, Warner's proposal is actually
optimal in a strong sense.  There are also other relatively simple
estimation strategies that achieve convergence rate $d
/ n \diffp^2$; for instance, as we show, the Laplace perturbation 
approach proposed by \citet{DworkMcNiSm06} is one such strategy. Nonetheless,
its ease of use coupled with our optimality results provide
support for randomized response as a preferred method 
for private estimation of population probabilities.

Let us now prove that these strategies attain the optimal rate
of convergence.  Since there is a bijection between multinomial
samples $\statsample \in \{1, \ldots, d\}$ and the $d$ standard basis
vectors $e_1, \ldots, e_d \in \R^d$, we abuse notation and represent
samples $\statsample$ as either when designing estimation strategies.

\paragraph{Randomized response:}

In randomized response, we construct the private vector $\channelrv
\in \{0, 1\}^d$ from a multinomial sample $\statsample \in \{e_1,
\ldots, e_d\}$ by sampling $d$ coordinates independently via the
procedure
\begin{equation}
  [\channelrv]_j = \begin{cases} \statsample_j & \mbox{with~probability~}
    \frac{\exp(\diffp / 2)}{1 + \exp(\diffp / 2)} \\
    1 - \statsample_j & \mbox{with~probability~}
    \frac{1}{1 + \exp(\diffp/2)}. \end{cases}
  \label{eqn:randomized-multinomial-response}
\end{equation}
We claim that this channel~\eqref{eqn:randomized-multinomial-response}
is $\diffp$-differentially private: indeed, note that for any
$\statsample, \statsample' \in \simplex_d$ and any vector $\channelval
\in \{0, 1\}^d$ we have
\begin{align*}
  \frac{\channelprob( \channelrv = \channelval \mid \statsample)}{
    \channelprob( \channelrv = \channelval \mid \statsample')} & =
  \frac{ \left( \frac{\exp(\diffp / 2)}{1 + \exp(\diffp / 2)}
    \right)^{ \lone{\channelval - \statsample}} \left( \frac{1}{1 +
      \exp(\diffp / 2)}\right)^{ d - \lone{\channelval -
        \statsample}}}{ \left(\frac{\exp(\diffp / 2)}{ 1 + \exp(\diffp
      / 2)}\right)^{ \lone{\channelval - \statsample'}} \left(
    \frac{1}{1 + \exp(\diffp / 2)}\right)^{ d - \lone{ \channelval -
        \statsample'}}} \\
& = \exp\left(\frac{\diffp}{2} \left(\lone{\channelval - \statsample}
  - \lone{\channelval - \statsample'}\right)\right) \in
  \left[\exp(-\diffp), \exp(\diffp)\right],
\end{align*}
where we used the triangle inequality to assert that
$|\lone{\channelval - \statsample} - \lone{\channelval -
  \statsample'}| \le \lone{\statsample - \statsample'} \le 2$.  We can
compute the expected value and variance of the random variables
$\channelrv$; indeed, by
definition~\eqref{eqn:randomized-multinomial-response}
\begin{equation*}
  \E[\channelrv \mid \statsample]
  = \frac{e^{\diffp / 2}}{1 + e^{\diffp/2}}
  \statsample + \frac{1}{1 + e^{\diffp/2}}
  (\onevec - \statsample)
  = \frac{e^{\diffp/2} - 1}{e^{\diffp/2} + 1}
  \statsample + \frac{1}{1 + e^{\diffp/2}} \onevec.
\end{equation*}
Since the $\channelrv$ are Bernoulli, we obtain
the variance bound $\E[\ltwo{\channelrv - \E[\channelrv]}^2]
< d/4 + 1$.
Recalling the definition of the projection
$\project_{\simplex_d}$ onto the simplex, we arrive at the natural estimator
\begin{equation}
  \what{\optvar}_{\rm part}
  \defeq \frac{1}{n}
  \sum_{i=1}^n \left(\channelrv_i - \onevec / (1 + e^{\diffp/2}) \right)
  \frac{e^{\diffp/2} + 1}{e^{\diffp/2} - 1}
  ~~~ \mbox{and} ~~~
  \what{\optvar} \defeq \project_{\simplex_d}\left(\what{\optvar}_{\rm part}
  \right).
  \label{eqn:randomized-response-estimator}
\end{equation}
The projection of $\what{\optvar}_{\rm part}$ onto the probability simplex can
be done in time linear in the dimension $d$ of the problem~\cite{Brucker84},
so the estimator~\eqref{eqn:randomized-response-estimator}
is efficiently computable.  Since projections only decrease distance, vectors
in the simplex are at most distance $\sqrt{2}$ apart, and
$\E_\optvar[\what{\optvar}_{\rm part}] = \optvar$ by construction, we find
that
\begin{align}
  \E \left [\ltwos{\what{\optvar} - \optvar}^2\right] \leq \min
  \left\{2, \E \left[ \ltwos{\what{\optvar}_{\rm part} - \optvar}^2
    \right] \right \} \le \min \bigg\{2, \left(\frac{d}{4n} +
  \frac{1}{n} \right ) \left ( \frac{e^{\diffp/2} + 1}{ e^{\diffp/2} -
    1} \right)^2\bigg\} \lesssim \min\left\{1, \frac{d}{n \diffp^2}
  \right \}.  \nonumber
\end{align}


\paragraph{Laplace perturbation:}
We now turn to the strategy of~\citet{DworkMcNiSm06}, where we add Laplacian 
noise to the data. Let the vector $W \in \R^d$ have independent coordinates, 
each distributed as $\laplace(\diffp / 2)$. Then $\statsample + W$ is
$\diffp$-differentially private for $\statsample \in \simplex_d$: the
ratio of the densities $\channeldensity(\cdot \mid \statsample)$ and
$\channeldensity(\cdot \mid \statsample')$ of $\statsample + W$ and
$\statsample' + W$ is
\begin{equation*}
  \frac{\channeldensity(\channelval \mid \statsample)}{
    \channeldensity(\channelval \mid \statsample')} =
  \frac{\exp(-(\diffp/2) \lone{\statsample - \channelval})}{
    \exp(-(\diffp/2) \lone{\statsample' - \channelval})} \in
  \left[\exp(-\diffp), \exp(\diffp)\right].
\end{equation*}
By defining the private data $\channelrv_i = \statrv_i + W_i$,
where $W_i \in \R^d$ are independent, we can define the partial
estimator $\what{\optvar}_{\rm part} = \frac{1}{n} \sum_{i=1}^n
\channelrv_i$ and the projected estimator $\what{\optvar} \defeq
\project_{\simplex_d}(\what{\optvar}_{\rm part})$, similar to our
randomized response
construction~\eqref{eqn:randomized-response-estimator}.  Then by
computing the variance of the noise samples $W_i$, it is clear that
\begin{equation}
  \E\left[\ltwos{\what{\optvar} - \optvar}^2\right] \le \min\left\{2,
  \E\left[\ltwos{\what{\optvar}_{\rm part} - \optvar}^2
    \right]\right\} \le \min\left\{2, \frac{1}{n} + 4 \frac{d}{n
    \diffp^2}\right\} \lesssim \min\left\{1, \frac{d}{n
    \diffp^2}\right\}.  \nonumber
\end{equation}
For small $\diffp$, the Laplace perturbation approach has a somewhat
sharper convergence rate in terms of constants than that of the
randomized response
estimator~\eqref{eqn:randomized-response-estimator}, so in some cases
it may be preferred. Nonetheless, the simplicity of explaining the
sampling procedure~\eqref{eqn:randomized-multinomial-response} may
argue for its use in scenarios such as survey sampling.


\section{Density Estimation under Local Privacy}
\label{sec:density-estimation}

In this section, we turn to studying a nonparametric statistical
problem in which the effects of local differential privacy turn out
to be somewhat more severe.  We show that for the problem of density
estimation, instead of just multiplicative loss in the effective
sample size as in previous section (see also our
paper~\cite{DuchiJoWa13_parametric}), imposing local differential
privacy leads to a completely different convergence rate.  This result
holds even though we solve an estimation problem in which the function
estimated and the samples themselves belong to compact spaces.

In more detail, we consider estimation of probability densities $f :
\R \rightarrow \R_+$, $\int f(x) dx = 1$ and $f \ge 0$, defined on the
real line, focusing on a standard family of densities of varying
smoothness~\cite[e.g.][]{Tsybakov09}.  Throughout this section, we let
$\numderiv \in \N$ denote a fixed positive integer.  Roughly, we
consider densities that have bounded $\numderiv$th derivative, and
we study density estimation using the squared $L^2$-norm $\ltwo{f}^2
\defeq \int f^2(x) dx$ as our metric; in formal terms, we impose these
constraints in terms of Sobolev classes
(e.g.~\cite{Tsybakov09,Efromovich99}). Let the countable collection of
functions $\{\basisfunc_j\}_{j=1}^\infty$ be an orthonormal basis for
$L^2([0, 1])$.  Then any function $f \in L^2([0, 1])$ can be expanded
as a sum $\sum_{j=1}^\infty \optvar_j \basisfunc_j$ in terms of the
basis coefficients $\optvar_j \defeq \int f(x) \basisfunc_j(x) dx$. By
Parseval's theorem, we are guaranteed that $\{\optvar_j\}_{j=1}^\infty
\in \ell^2(\N)$.  The \emph{Sobolev space} $\densclass[\lipconst]$ is
obtained by enforcing a particular decay rate on the coefficients
$\optvar$:
\begin{definition}[Elliptical Sobolev space]
\label{definition:sobolev-densities}
For a given orthonormal basis $\{\basisfunc_j\}$ of $L^2([0, 1])$,
smoothness parameter $\beta > 1/2$ and radius $\lipconst$, the
function class $\densclass[\lipconst]$ is given by
  \begin{equation*}
    \densclass[\lipconst] \defeq \bigg\{f \in L^2([0, 1]) \mid f =
    \sum_{j=1}^\infty \optvar_j \basisfunc_j ~ \mbox{such~that} ~
    \sum_{j=1}^\infty j^{2\numderiv} \basisfunc_j^2 \le
    \lipconst^2\bigg\}.
  \end{equation*}
\end{definition}

If we choose the trignometric basis as our orthonormal basis, then
membership in the class $\densclass[\lipconst]$ corresponds to certain
smoothness constraints on the derivatives of $f$.  More precisely, for
$j \in \N$, consider the orthonormal basis for $L^2([0, 1])$ of
trigonometric functions:
\begin{equation}
  \label{eqn:trig-basis}
  \basisfunc_0(t) = 1,
  ~~~ \basisfunc_{2j}(t) = \sqrt{2} \cos(2\pi j t),
  ~~~ \basisfunc_{2j + 1}(t) = \sqrt{2} \sin(2\pi j t).
\end{equation}
Now consider a $\numderiv$-times almost everywhere differentiable
function $f$ for which $|f^{(\numderiv)}(x)| \le \lipconst$ for almost
every $x \in [0, 1]$ satisfying $f^{(k)}(0) = f^{(k)}(1)$ for $k \le
\numderiv - 1$. Then it is known~\cite[Lemma A.3]{Tsybakov09} that,
uniformly for such $f$, there is a universal constant $c$ such that
that $f \in \densclass[c \lipconst]$.
Thus, Definition~\ref{definition:sobolev-densities} (essentially)
captures densities that have Lipschitz-continuous $(\numderiv - 1)$th
derivative. In the sequel, we write $\densclass$ when the bound
$\lipconst$ in $\densclass[\lipconst]$ is $\order(1)$. It is well
known~\cite{Yu97,YangBa99,Tsybakov09} that the minimax risk for
non-private estimation of densities in the class $\densclass$ scales
as
\begin{equation}
\label{eqn:classical-density-estimation-rate}
\minimax_n \left( \densclass, \ltwo{\cdot}^2, \infty\right) \asymp
n^{-\frac{2 \numderiv}{2 \numderiv + 1}}.
\end{equation}
The goal of this section is to understand how this minimax rate
changes when we add an $\diffp$-privacy constraint to the problem.  Our
main result is to demonstrate that the classical
rate~\eqref{eqn:classical-density-estimation-rate} is \emph{no longer
  attainable} when we require $\diffp$-local differential privacy.  In
particular, we prove a lower bound that is substantially larger.
In Sections~\ref{sec:histogram-estimators}
and~\ref{sec:orthogonal-series}, we show how to achieve this lower
bound using histogram and orthogonal series estimators.

\subsection{Lower bounds on density estimation}

We begin by giving our main lower bound on the minimax rate of estimation of
densities when samples from the density must be kept differentially private.
We provide the proof of the following theorem in
Section~\ref{sec:proof-density-estimation}.
\begin{theorem}
  \label{theorem:density-estimation}
  Consider the class of densities $\densclass$ defined using the
  trigonometric basis~\eqref{eqn:trig-basis}.  For some $\diffp \in
  (0, 1/4]$, suppose $\channelrv_i$ are $\diffp$-locally
    private~\eqref{eqn:local-privacy} for the samples $\statrv_i \in
    [0, 1]$.  There exists a constant $c > 0$, dependent only on
    $\numderiv$, such that
\begin{align}
 \label{eqn:private-density-estimation-rate}
 \minimax_n \left( \densclass, \ltwo{\cdot}^2, \diffp \right) & \geq c
 \left (n \diffp^2 \right)^{-\frac{2 \numderiv}{2 \numderiv + 2}}.
\end{align}
\end{theorem}
\noindent In comparison with the classical minimax
rate~\eqref{eqn:classical-density-estimation-rate}, the lower
bound~\eqref{eqn:private-density-estimation-rate} is substantially
different, in that it involves a different polynomial exponent:
namely, the exponent is $2 \numderiv / (2 \numderiv + 1)$ in the
classical case~\eqref{eqn:classical-density-estimation-rate}, while in
the differentially private
case~\eqref{eqn:private-density-estimation-rate}, the exponent has
been reduced to $2 \numderiv / (2 \numderiv + 2)$.  For example, when
we estimate Lipschitz densities, we have $\numderiv = 1$, and the rate
degrades from $n^{-2/3}$ to $n^{-1/2}$.  Moreover, this degradation
occurs even though our samples are drawn from a compact space and the
set $\densclass$ is also compact.  

Interestingly, no estimator based on Laplace (or exponential)
perturbation of the samples $\statrv_i$ themselves can attain the rate
of convergence~\eqref{eqn:private-density-estimation-rate}.  This
can be established by connecting such a perturbation-based approach to
the problem of nonparametric deconvolution.  In their study of the
deconvolution problem, \citet{CarrollHa88} show that if samples
$\statrv_i$ are perturbed by additive noise $W$, where the
characteristic function $\phi_W$ of the additive noise has tails
behaving as $|\phi_W(t)| = \order(|t|^{-a})$ for some $a > 0$, then no
estimator can deconvolve the samples $\statrv + W$ and attain a rate
of convergence better than
$n^{-2 \numderiv / (2 \numderiv + 2a + 1)}$.
Since the Laplace distribution's characteristic function has tails
decaying as $t^{-2}$, no estimator based on perturbing the samples
directly can attain a rate of convergence better than $n^{-2\numderiv
  / (2 \numderiv + 5)}$.  If the lower
bound~\eqref{eqn:private-density-estimation-rate} is attainable, we
must then study privacy mechanisms that are not simply based on direct
perturbation of the samples $\{\statrv_i\}_{i=1}^\numobs$.


\subsection{Achievability by histogram estimators}
\label{sec:histogram-estimators}

We now turn to the mean-squared errors achieved by specific practical
schemes, beginning with the special case of Lipschitz density
functions ($\numderiv = 1$), for which it suffices to consider a
private version of a classical histogram estimate.

For a fixed positive integer $\numbin \in \N$, let
$\{\statdomain_j\}_{j=1}^\numbin$ denote the partition of $\statdomain
= [0, 1]$ into the intervals 
\begin{align*}
  \statdomain_j = \openright{(j - 1) / \numbin}{j/\numbin} \quad
  \mbox{for $j = 1, 2, \ldots, \numbin-1$,~~and} ~~
  \statdomain_\numbin = [(\numbin-1)/\numbin, 1].
\end{align*}
Any histogram estimate of the density based on these $\numbin$ bins
can be specified by a vector $\optvar \in \numbin
\simplex_\numbin$, where we recall $\simplex_\numbin \subset
\R^\numbin_+$ is the probability simplex. Any such vector 
defines a density estimate via the sum
\begin{equation*}
  f_\optvar \defeq \sum_{j=1}^\numbin \optvar_j
  \characteristic{\statdomain_j},
\end{equation*}
where $\characteristic{E}$ denotes the characteristic (indicator) function of
the set $E$.  

Let us now describe a mechanism that guarantees $\diffp$-local
differential privacy.  Given a data set $\{\statrv_1, \ldots,
\statrv_n\}$ of samples from the distribution $f$, consider the
vectors
\begin{align}
\channelrv_i & \defeq \histelement_\numbin(\statrv_i) + W_i, \quad
\mbox{for $i = 1, 2, \ldots, \numobs$},
\end{align}
where $\histelement_\numbin(\statrv_i) \in \simplex_\numbin$ is a
$\numbin$-vector with the $j^{th}$ entry equal to one if $\statrv_i
\in \statdomain_j$, and zeroes in all other entries, and $W_i$ is a
random vector with i.i.d.\ $\laplace(\diffp/2)$ entries.
The variables $\{\channelrv_i\}_{i=1}^\numobs$ so-defined
are $\diffp$-locally differentially private for
$\{\statrv_i\}_{i=1}^\numobs$.

Using these private variables, we then form the density estimate
$\what{f} \defeq f_{\what{\optvar}} = \sum_{j=1}^\numbin
\what{\optvar}_j \characteristic{\statdomain_j}$ based on the vector
\begin{align}
  \label{eqn:histogram-estimator}
  \what{\optvar} & \defeq \Pi_\numbin \bigg(\frac{\numbin}{n}
  \sum_{i=1}^n \channelrv_i\bigg),
\end{align}
where $\Pi_\numbin$ denotes the Euclidean projection operator onto the
set $\numbin \simplex_\numbin$.  By construction, we have $\what{f}
\ge 0$ and $\int_0^1 \what{f}(\statsample) d\statsample = 1$, so
$\what{f}$ is a valid density estimate.
\begin{proposition}
  \label{proposition:histogram-estimator}
  Consider the estimate $\what{f}$ based on $\numbin = (n \diffp^2)^{1/4}$
  bins in the histogram.  For any $1$-Lipschitz density $f : [0,
    1] \rightarrow \R_+$, we have
  \begin{align}
    \label{EqnHistoAchieve}
    \E_f\left [ \ltwobig{\what{f} - f}^2 \right ] & \le 5 (\diffp^2
    n)^{-\half} + \sqrt{\diffp} n^{-3/4}.
  \end{align}
\end{proposition}
\noindent For any fixed $\diffp > 0$, the first term in the
bound~\eqref{EqnHistoAchieve} dominates, and the $\order
((\diffp^2 \numobs)^{-\half})$ rate matches the minimax
lower bound~\eqref{eqn:private-density-estimation-rate} in the case
$\beta = 1$.  Consequently, we have shown that the privatized
histogram estimator is minimax-optimal for Lipschitz densities.  This
result provides the private analog of the classical result that
histogram estimators are minimax-optimal (in the non-private setting)
for Lipshitz densities.  See Section~\ref{sec:proof-histogram} for a
proof of Proposition~\ref{proposition:histogram-estimator}. We
remark in passing that a randomized response scheme parallel to
that of Section~\ref{sec:private-multinomial-estimation} achieves
the same rate of convergence; once again, randomized response is optimal.


\subsection{Achievability by orthogonal projection estimators}
\label{sec:orthogonal-series}

For higher degrees of smoothness ($\numderiv > 1$), histogram
estimators no longer achieve optimal rates in the classical setting.
Accordingly, we now turn to developing estimators based on orthogonal
series expansion, and show that even in the setting of local privacy,
they can achieve the lower
bound~\eqref{eqn:private-density-estimation-rate} for all orders of
smoothness $\numderiv \ge 1$.  

Recall the elliptical Sobolev space
(Definition~\ref{definition:sobolev-densities}), in which a function
$f$ is represented in terms of its basis expansion $f = \sum_{j =
  1}^\infty \optvar_j \basisfunc_j$, where $\optvar_j = \int f(x)
\basisfunc_j(x) dx$.  This representation underlies the classical
method of orthonormal series estimation: given a data set, 
$\{\statrv_1, \statrv_2, \ldots, \statrv_n\}$, drawn i.i.d.\ according 
to a density $f \in L^2([0, 1])$, we first compute the
empirical basis coefficients
\begin{equation}
\what{\optvar}_j = \frac{1}{n} \sum_{i=1}^n \basisfunc_j(\statrv_i)
~~~ \mbox{and then set} ~~~ \what{f} = \sum_{j=1}^\numbin
\what{\optvar}_j \basisfunc_j,
  \label{eqn:projection-estimator}
\end{equation}
where the value $\numbin \in \N$ is chosen either a priori based on
known properties of the estimation problem or adaptively, for example,
using cross-validation~\cite{Efromovich99,Tsybakov09}.

In the setting of local privacy, we consider a mechanism that, instead of
releasing the vector of coefficients $\big(\basisfunc_1(\statrv_i), \ldots,
\basisfunc_\numbin(\statrv_i) \big)$ for each data point, employs a random vector
$\channelrv_i = (\channelrv_{i,1}, \ldots, \channelrv_{i, \numbin})$ with the
property that $\E [ \channelrv_{i,j} \mid \statrv_i] =
\basisfunc_j(\statrv_i)$ for each $j = 1, 2, \ldots, \numbin$.  We assume the
basis functions are uniformly bounded; i.e., there exists a constant
$\orthbound = \sup_j \sup_\statsample |\basisfunc_j(\statsample)| < \infty$.
(This boundedness condition holds for many standard bases, including the
trigonometric basis underlying the classical Sobolev classes and the Walsh
basis.) For a fixed number $\sbound$ strictly larger than $\orthbound$ (to
be specified momentarily), consider the following scheme:
\begin{description}
\item[Sampling strategy] Given a vector $\tau \in [-\orthbound,
  \orthbound]^\numbin$, construct $\wt{\tau} \in \{-\orthbound,
  \orthbound\}^\numbin$ with coordinates $\wt{\tau}_j$ sampled
  independently from $\{-\orthbound, \orthbound\}$ with probabilities
  $\half - \frac{\tau_j}{2 \orthbound}$ and $\half + \frac{\tau_j}{2
    \orthbound}$.  Sample $T$ from a $\bernoulli(e^\diffp / (e^\diffp
  + 1))$ distribution. Then choose $\channelrv \in \{-\sbound,
  \sbound\}^\numbin$ via
\begin{equation}
 \label{eqn:infty-halfspace-samples}
\channelrv \sim \begin{cases} \mbox{Uniform~on~} \left \{\channelval
  \in \{-\sbound, \sbound\}^\numbin : \<\channelval, \wt{\tau}\> > 0
  \right \} & \mbox{if~} T = 1 \\
\mbox{Uniform~on~} \left \{ \channelval \in \{-\sbound,
\sbound\}^\numbin : \< \channelval, \wt{\tau} \> \le 0 \right\} &
\mbox{if~} T = 0.
    \end{cases}
  \end{equation}
\end{description}
By inspection, $\channelrv$ is $\diffp$-differentially private for any
initial vector in the box $[-\orthbound, \orthbound]^\numbin$, and
moreover, the samples~\eqref{eqn:infty-halfspace-samples} are
efficiently computable (for example by rejection sampling).

Starting from the vector $\tau \in \R^\numbin$, $\tau_j =
\basisfunc_j(\statrv_i)$, in the above sampling strategy, iteration of
expectation yields
\begin{equation}
  \E[[\channelrv]_j \mid \statrv = \statsample]
  = c_\numbin \frac{\sbound}{\orthbound \sqrt{\numbin}}
  \left(\frac{e^\diffp}{e^\diffp + 1} - \frac{1}{e^\diffp + 1}\right)
  \basisfunc_j(\statsample)
  = c_\numbin \frac{\sbound}{\orthbound \sqrt{\numbin}} \frac{e^\diffp - 1}{
    e^\diffp + 1} \basisfunc_j(\statsample), 
  \label{eqn:size-infinity-channel}
\end{equation}
for a constant $c_\numbin$ that may depend on $\numbin$ but is $\order(1)$ and
bounded away from 0.  Consequently, to attain the unbiasedness condition
$\E[[\channelrv_i]_j \mid \statrv_i] = \basisfunc_j(\statrv_i)$, it suffices
to take $\sbound = \order(\orthbound \sqrt{\numbin} / \diffp)$.

The full sampling and inferential scheme are as follows: given a data point
$\statrv_i$, we sample $\channelrv_i$ according to the
strategy~\eqref{eqn:infty-halfspace-samples}, where we start from the vector
$\tau = [\basisfunc_j(\statrv_i)]_{j=1}^\numbin$ and use the bound $\sbound =
\orthbound \sqrt{\numbin} (e^\diffp + 1) / c_\numbin (e^\diffp - 1)$, where
the constant $c_\numbin$ is as in the
expression~\eqref{eqn:size-infinity-channel}.  Defining the density estimator
\begin{equation}
  \what{f} \defeq \frac{1}{\numobs} \sum_{i=1}^\numobs \sum_{j=1}^\numbin
  \channelrv_{i,j} \basisfunc_j,
  \label{eqn:orthogonal-density-estimator}
\end{equation}
we obtain the following proposition.
\begin{proposition}
  \label{proposition:density-upper-bound}
  Let $\{\basisfunc_j\}$ be a $\orthbound$-uniformly bounded
  orthonormal basis for $L^2([0, 1])$.  There exists a constant
  $c$ (depending only on $\lipconst$ and $\orthbound$) such that the
  estimator~\eqref{eqn:orthogonal-density-estimator} with
  \mbox{$\numbin = (n \diffp^2)^{1 / (2 \numderiv + 2)}$} satisfies
  \begin{equation*}
    \E_f \left [\ltwos{f - \what{f}}^2 \right] \leq c \left (n \diffp^2
    \right)^{-\frac{2 \numderiv}{2 \numderiv + 2}}.
  \end{equation*}
  for any $f$ in the Sobolev space $\densclass[\lipconst]$.
\end{proposition}
\noindent
See Section~\ref{sec:proof-density-upper-bound} for a proof.

Propositions~\ref{proposition:histogram-estimator} and
\ref{proposition:density-upper-bound} make clear that the minimax lower
bound~\eqref{eqn:private-density-estimation-rate} is sharp, as claimed. We
have thus attained a variant of the known minimax density estimation
rates~\eqref{eqn:classical-density-estimation-rate}, but with a polynomially
worse sample complexity as soon as we require local differential privacy.

Before concluding our exposition, we make a few remarks on other potential
density estimators.  Our orthogonal-series
estimator~\eqref{eqn:orthogonal-density-estimator} (and sampling
scheme~\eqref{eqn:size-infinity-channel}), while similar in spirit to that
proposed by~\citet[Sec.~6]{WassermanZh10}, is different in that it is locally
private and requires a different noise strategy to obtain both $\diffp$-local
privacy and optimal convergence rate.  Finally, we consider the insufficiency
of standard Laplace noise addition for estimation in the setting of this
section. Consider the vector \mbox{$[\basisfunc_j(\statrv_i)]_{j=1}^\numbin
  \in [-\orthbound, \orthbound]^\numbin$.} To make this vector
$\diffp$-differentially private by adding an independent Laplace noise vector
$W \in \R^\numbin$, we must take $W_j \sim \laplace(\diffp / (\orthbound
k))$. The natural orthogonal series estimator (e.g.~\cite{WassermanZh10}) is
to take $\channelrv_i = [\basisfunc_j(\statrv_i)]_{j=1}^\numbin + W_i$, where
$W_i \in \R^\numbin$ are independent Laplace noise vectors. We then use the
density estimator~\eqref{eqn:orthogonal-density-estimator}, except that we use
the Laplacian perturbed $\channelrv_i$. However, this estimator suffers the
following drawback (see section~\ref{sec:proof-laplace-density-bad}):
\begin{observation}
  \label{observation:laplace-density-bad}
  Let $\what{f} = \frac{1}{\numobs} \sum_{i=1}^\numobs \sum_{j=1}^\numbin
  \channelrv_{i, j} \basisfunc_j$, where the $\channelrv_i$ are the
  Laplace-perturbed vectors of the previous paragraph. Assume the orthonormal
  basis $\{\basisfunc_j\}$ of $L^2([0, 1])$ contains the constant function.
  There is a constant $c$ such that for any $\numbin \in \N$,
  there is an $f \in \densclass[2]$ such that
  \begin{equation*}
    \E_f\left[\ltwos{f - \what{f}}^2\right] \ge 
    c (n \diffp^2)^{-\frac{2 \numderiv}{2 \numderiv + 3}}.
  \end{equation*}
\end{observation}
\noindent
This lower bound shows that standard estimators
based on adding Laplace noise to appropriate basis expansions of the
data fail: there is a degradation in rate from
$n^{-\frac{2 \numderiv}{2 \numderiv + 2}}$ to $n^{-\frac{2
    \numderiv}{2 \numderiv + 3}}$. While this is not a formal proof
that no approach based on Laplace perturbation can provide optimal
convergence rates in our setting, it does suggest that finding such an
estimator is non-trivial.


\section{Proof of Theorem~\ref{theorem:multinomial-rate}}
\label{sec:proof-multinomial-rate}

At a high level, our proof can be split into three steps, the first of
which is relatively standard, while the second two exploit specific
aspects of the local privacy set-up:

\begin{enumerate}[(1)]
\item 
\label{step:standard}
The first step is a standard reduction, based on
  Lemma~\ref{lemma:minimax-risk-bound}, from an estimation problem to
  a multi-way testing problem that involves discriminating between
  indices $\packval$ contained within some subset $\packset$ of
  $\R^d$.
\item \label{item:construct-packing} The second step is an appropriate
  construction of the set $\packset \subset \R^d$ such that each pair
  is $\delta$-separated and the resulting set is as large as possible
  (a maximal $\delta$-packing).  In addition, our arguments require
  that, for a random variable $\packrv$ uniformly distributed over
  $\packset$, the covariance $\cov(\packrv)$ has relatively small
  operator norm.
\item The final step is to apply
  Proposition~\ref{proposition:information-bounds} in order to control the
  mutual information associated with the testing problem. To do so, it is
  necessary to show that controlling the supremum subsets of
  $L^\infty(\statdomain)$ in the bound~\eqref{eqn:super-fano} can be reduced
  to bounding the operator norm of $\cov(\packrv)$.
\end{enumerate}

\noindent We have already described the reduction of
Step~\ref{step:standard} in Section~\ref{sec:minimax-framework}.
Accordingly, we turn to the second step.

\paragraph{Constructing a good packing:}

The following result on the binary hypercube \mbox{$H_d
  \defeq \{0, 1\}^d$} underlies our construction:
\begin{lemma}
  \label{lemma:lots-of-packings}
  There exist universal constants $c_1, c_2 \in (0, \infty)$ such that
  for each $\nnz \in \{1, \ldots, d\}$, there is a set $\packset \subset
  H_d$ with the following properties:
  \begin{enumerate}[(i)]
  \item \label{item:nnz}
    Any $\packval \in \packset$ has exactly $\nnz$ non-zero
    entries.
  \item \label{item:separate}
    For any $\packval, \altpackval \in \packset$ with $\packval
    \neq \altpackval$, the $\ell_1$-norm is lower bounded as
    $\lone{\packval - \altpackval} \ge \max\{\floor{\nnz/4}, 1\}$.
  \item The set $\packset$ has cardinality $\card(\packset) \ge
    (d/\nnz)^{c_1 \nnz}$.
  \item \label{item:sparse-covariance}
    For a random vector $\packrv$ uniformly distributed over
    $\packset$, we have
    \begin{align*}
      \cov(\packrv) & \preceq c_2 \frac{\nnz}{d} I_{d \times d}.
    \end{align*}
  \end{enumerate}
\end{lemma}
\noindent The proof of Lemma~\ref{lemma:lots-of-packings} is based on the
probabilistic method~\cite{AlonSp00}: we show that a certain randomized
procedure generates such a packing with strictly positive probability.  Along
the way, we use matrix Bernstein inequalities~\cite{MackeyJoChFaTr12} and some
approximation-theoretic ideas developed by~\citet{Kuhn01}. We provide details
in Appendix~\ref{appendix:proof-lemma-lots-of-packings}. \\

We now construct a suitable
packing of the the unit simplex $\simplex_d$.  Given an integer $\nnz
\in \{1, \ldots, d\}$, consider the packing $\packset \subset \{0,
1\}^d$ given by Lemma~\ref{lemma:lots-of-packings}.  For a fixed
$\delta \in [0, 1]$, consider the following family of vectors in $\R^d$
\begin{equation*}
  \optvar_\packval \defeq \frac{\delta}{\nnz} \packval + \frac{(1 -
    \delta)}{d} \onevec, \quad \mbox{for each $\packval \in
    \packset$.}
\end{equation*}
By inspection, each of these vectors belongs to the $d$-variate probability
simplex (i.e., satisfies $\<\onevec, \optvar_\packval\> = 1$ and
$\optvar_\packval \geq 0$). Moreover, since the vector $\packval -
\altpackval$ can have at most $2 \nnz$ non-zero entries, we have $\|\packval -
\altpackval\|_1 \leq \sqrt{2 \nnz} \|\packval - \altpackval\|_2$.  Combined
with property~(\ref{item:separate}), we conclude that
for universal constants $c, c' > 0$
\begin{align*}
  \ltwo{\packval - \altpackval}
  \ge \frac{\lone{\packval - \altpackval}}{\sqrt{2 \nnz}}
  \geq c' \nnz \frac{1}{\sqrt{2
    \nnz}} = c \sqrt{\nnz}.
\end{align*}
By the definition of $\optvar_\packval$, we then have for a universal
constant $c$ that
\begin{align}
  \label{EqnKeyLower}
  \ltwo{\optvar_\packval - \optvar_{\altpackval}}^2 =
  \frac{\delta^2}{\nnz^2} \ltwo{\packval - \altpackval}^2 & \geq
  c \frac{\delta^2}{\nnz}.
\end{align}


\paragraph{Upper bounding the mutual information:}

Our next step is to upper bound the mutual information
$\information(\channelrv_1, \ldots, \channelrv_n; \packrv)$.  Recall
the definition of the linear functionals $\varphi_\packval$ from
Proposition~\ref{proposition:information-bounds}.  Since $\statdomain
= \{1, 2, \ldots, d\}$, any element of $L^\infty(\statdomain)$ may be
identified with a vector $\optdens \in \R^d$.  Following this
identification, we have
\begin{align*}
\varphi_\packval(\optdens) = \sum_{j=1}^d \optvar_{\packval, j}
\optdens_j - \frac{1}{|\packset|} \sum_{\altpackval \in \packset}
\sum_{j=1}^d \optvar_{\altpackval, j} \optdens_j = \frac{\delta}{\nnz}
\<\optdens, \packval - \E[\packrv]\>,
\end{align*}
where $\packrv$ is a random variable distributed uniformly over
$\packset$. As a consequence, we have
\begin{equation*}
  \frac{1}{|\packset|} \sum_{\packval \in \packset}
  \left(\varphi_\packval(\optdens)\right)^2
  = \frac{\delta^2}{\nnz^2} \optdens^\top \cov(V) \optdens
  \le c_2 \frac{\delta^2}{d\nnz} \ltwo{\optdens}^2,
\end{equation*}
where the final inequality follows from
Lemma~\ref{lemma:lots-of-packings}(\ref{item:sparse-covariance}). For any
$\optdens \in \linfset$, we have the upper bound $\ltwo{\optdens}^2 \le d
(e^\diffp - e^{-\diffp})^2 / 4$, whence
\begin{equation*}
\sup_{\optdens \in \linfset} \frac{1}{|\packset|} \sum_{\packval \in
  \packset} \left(\varphi_\packval(\optdens)\right)^2 \le c
\frac{(e^\diffp - e^{-\diffp})^2 \delta^2}{\nnz},
\end{equation*}
for some universal constant $c$.  Consequently, by applying the
information inequality~\eqref{eqn:super-fano}, we conclude
that there is a universal constant
constant $C$ such that
\begin{align}
  \information(\channelrv_1, \ldots, \channelrv_n; \packrv) & \leq C
  \frac{n \diffp^2 \delta^2}{\nnz} \qquad \mbox{for all $\diffp \le
    1/4$.}
\end{align}

\paragraph{Applying testing inequalities:}

The final step is to lower bound the testing error.  Since the vectors $\{
\optvar_\packval, \packval \in \packset\}$ are $c / \sqrt{\nnz}$-separated in
the $\ell_2$-norm~\eqref{EqnKeyLower} and Lemma~\ref{lemma:lots-of-packings}
implies $\card(\packset) \ge (d/\nnz)^{c_2 \nnz}$ for a constant $c_2$, Fano's
inequality~\eqref{eqn:fano} implies
\begin{equation}
  \label{eqn:multinomial-apply-fano}
  \minimax_n \left (\simplex_d, \ltwo{\cdot}^2, \diffp\right) \ge
  c_0 \frac{\delta^2}{\nnz} \left(1 - \frac{c_1 n \delta^2 \diffp^2 / \nnz
    + \log 2}{ c_2 \nnz \log (d / \nnz)}\right),
\end{equation}
for universal constants $c_0, c_1, c_2$.  We split the remainder of
our analysis into cases, depending on the values of $(\nnz, d)$.  \\

\noindent {\emph{Case 1:}} First, suppose that $(\nnz, d)$ are large
enough to guarantee that 
\begin{align}
  \label{EqnCaseSplit}
  c_2 \nnz \log (d/\nnz) \ge 3 \log 2.
\end{align}
In this case, if we set
\begin{equation*}
  \delta^2 = \min\left\{1, \frac{c_2 \nnz^2}{2 c_1 n \diffp^2}
  \log \frac{d}{\nnz} \right\},
\end{equation*}
then we have
\begin{align*}
  1 - \frac{c_1 n \delta^2 \diffp^2 / \nnz + \log 2}{c_2 \nnz \log (d /
    \nnz)} & \geq 1 - \frac{\frac{c_1 n \diffp^2}{\nnz} \frac{c_2
      \nnz^2}{2 c_1 n \diffp^2} \log \frac{d}{\nnz} + \log 2}{c_2 \nnz
    \log \frac{d}{\nnz}} \\
  & = 1 - \frac{\half c_2 \nnz \log \frac{d}{\nnz} + \log 2}{ c_2 \nnz
  \log \frac{d}{\nnz}} \ge 1 - \half - \frac{1}{3} = \frac{1}{6}.
\end{align*}
Combined with the Fano lower bound~\eqref{eqn:multinomial-apply-fano},
this yields the claim~\eqref{eqn:multinomial-rate} under
condition~\eqref{EqnCaseSplit}. \\

\noindent \emph{Case 2:} Alternatively, when $d$ is small enough that
condition~\eqref{EqnCaseSplit} is violated, we instead apply Le Cam's
inequality~\eqref{eqn:le-cam} to a two-point hypothesis.  For our
purposes, it suffices to consider the case $d = 2$, since for the
purposes of lower bounds, any higher-dimensional problem is at least
as hard as this case.  Define the two vectors
\begin{equation*}
  \optvar_1 = \frac{1 + \delta}{2} e_1 + \frac{1 - \delta}{2} e_2
  ~~~ \mbox{and} ~~~
  \optvar_2 = \frac{1 - \delta}{2} e_1 + \frac{1 + \delta}{2} e_2.
\end{equation*}
By construction, each of these vectors belongs to the probability simplex in
$\R^2$, and moreover, we have $\ltwo{\optvar_1 - \optvar_2}^2 = 2 \delta^2$.
Letting $\statprob_j$ denote the multinomial distribution defined by
$\optvar_j$, we also have $\tvnorm{\statprob_1 - \statprob_2} =
\lone{\optvar_1 - \optvar_2}/2 = \delta$.

In terms of the marginal measures $\marginprob_\packval^n$ defined in
equation~\eqref{eqn:marginal-channel}, Pinsker's inequality
(e.g.~\cite[Lemma 2.5]{Tsybakov09}) implies that
\begin{align*}
\tvnorm{\marginprob_1^n - \marginprob_2^n} \le
\sqrt{\dkl{\marginprob_1^n}{\marginprob_2^n} / 2} .
\end{align*}
Combined with Le Cam's inequality~\eqref{eqn:le-cam} and the upper
bound on KL divergences from
Proposition~\ref{proposition:information-bounds}, we find that the
minimax risk is lower bounded as
\begin{equation*}
\minimax_n(\simplex_d, \ltwo{\cdot}^2, \diffp) \ge \frac{\delta^2}{2}
\left(\half - \half \sqrt{2(e^\diffp - 1)^2 n \tvnorm{\statprob_1 -
    \statprob_2}^2} \right),
\end{equation*}
where $\statprob_i$ denotes the multinomial probability associated
with the vector $\optvar_i$.  Since $\tvnorm{\statprob_1 -
  \statprob_2} = \delta$ by construction and $e^\diffp - 1 \le (5/4)
\diffp$ for all $\diffp \in [0, 1/4]$, we have
\begin{equation*}
  \frac{\delta^2}{2} \left(\half - \half
  \sqrt{25 n \diffp^2 \delta^2 / 8}\right)
  = \frac{\delta^2}{2} \left(\half - \frac{5}{4 \sqrt{2}}
  \sqrt{n} \diffp \delta\right).
\end{equation*}
Choosing $\delta = \min\{1, \sqrt{2} / (5 \sqrt{n \diffp^2})\}$
guarantees that $\half - 5 \sqrt{n} \diffp \delta / 4 \sqrt{2} \ge
1/4$, and hence
\begin{equation*}
  \minimax_n(\simplex_d, \ltwo{\cdot}^2, \diffp) \ge \frac{1}{8}
  \min\left\{1, \frac{2}{25 n \diffp^2}\right\},
\end{equation*}
which completes the proof of Theorem~\ref{theorem:multinomial-rate}.

\paragraph{Proof of inequality~\eqref{eqn:multinomial-limits}}
We conclude by proving
inequality~\eqref{eqn:multinomial-limits}. We distinguish three cases:
\begin{equation*}
  \mbox{(i)} ~ n \diffp^2 < \log d,
  ~~ \mbox{(ii)} ~ \log d \le n \diffp^2
  \le \half d^2,
  ~~ \mbox{(iii)}~  n \diffp^2 \ge \half d^2.
\end{equation*}
In case (i), by taking $k = 1$ in the lower bound~\eqref{eqn:multinomial-rate}
we obtain the lower bound 1. In case (ii), we set
$k = \sqrt{n \diffp^2} \in [\sqrt{\log d}, d / \sqrt{2}]$, and we obtain
\begin{equation*}
  \min\left\{\frac{1}{k},
  \frac{k \log \frac{d}{k}}{n \diffp^2} \right\}
  = \min\left\{
  \frac{1}{\sqrt{n \diffp^2}},
  \frac{\log d}{\sqrt{n \diffp^2}}
  - \frac{\half \log(n \diffp^2)}{\sqrt{n \diffp^2}}\right\}
  \ge \frac{\log d}{\sqrt{n \diffp^2}}
  - \frac{\half (2 \log d + \log \half)}{\sqrt{n \diffp^2}}
  = \frac{\log 2}{2 \sqrt{n \diffp^2}}.
\end{equation*}
In the final case (iii), 
choosing $k = d / 2$ yields the bound
$\min\{2/d, d \log 2 / (n \diffp^2)\} \ge d \log 2 / (n \diffp^2)$.


\section{Proofs of Density Estimation Results}

In this section, we provide the proofs of the results stated in
Section~\ref{sec:density-estimation} on density estimation. We defer the
proofs of more technical results to the appendices.  Throughout all proofs, we
use $c$ to denote a universal constant whose value may change from line to
line.

\subsection{Proof of Theorem~\ref{theorem:density-estimation}}
\label{sec:proof-density-estimation}

As with our previous proof, the argument follows the general outline
described at the beginning of
Section~\ref{sec:proof-multinomial-rate}. We remark that our proof is
based on a local packing technique, a more classical approach than the
metric entropy approach developed by~\citet{YangBa99}.  We do so because
in the setting of local differential privacy we do not expect that global
results on metric entropy will be generally useful; rather, we must
carefully construct our packing set to control the mutual information,
relating the geometry of the packing to the actual information
communicated.  In comparison with our proof of
Theorem~\ref{theorem:multinomial-rate}, the construction of a suitable
packing of $\densclass$ is somewhat more challenging: the
identification of densities with finite-dimensional vectors,
which we require for our application of
Proposition~\ref{proposition:information-bounds}, is not
immediately obvious.  In all cases, we
use the trigonometric basis to prove our lower bounds, so
we may work directly with smooth density functions $f$.

\begin{figure}
  \begin{center}
    \begin{tabular}{cc}
      \psfrag{g}{$g_1$}
      \includegraphics[height=.32\columnwidth]{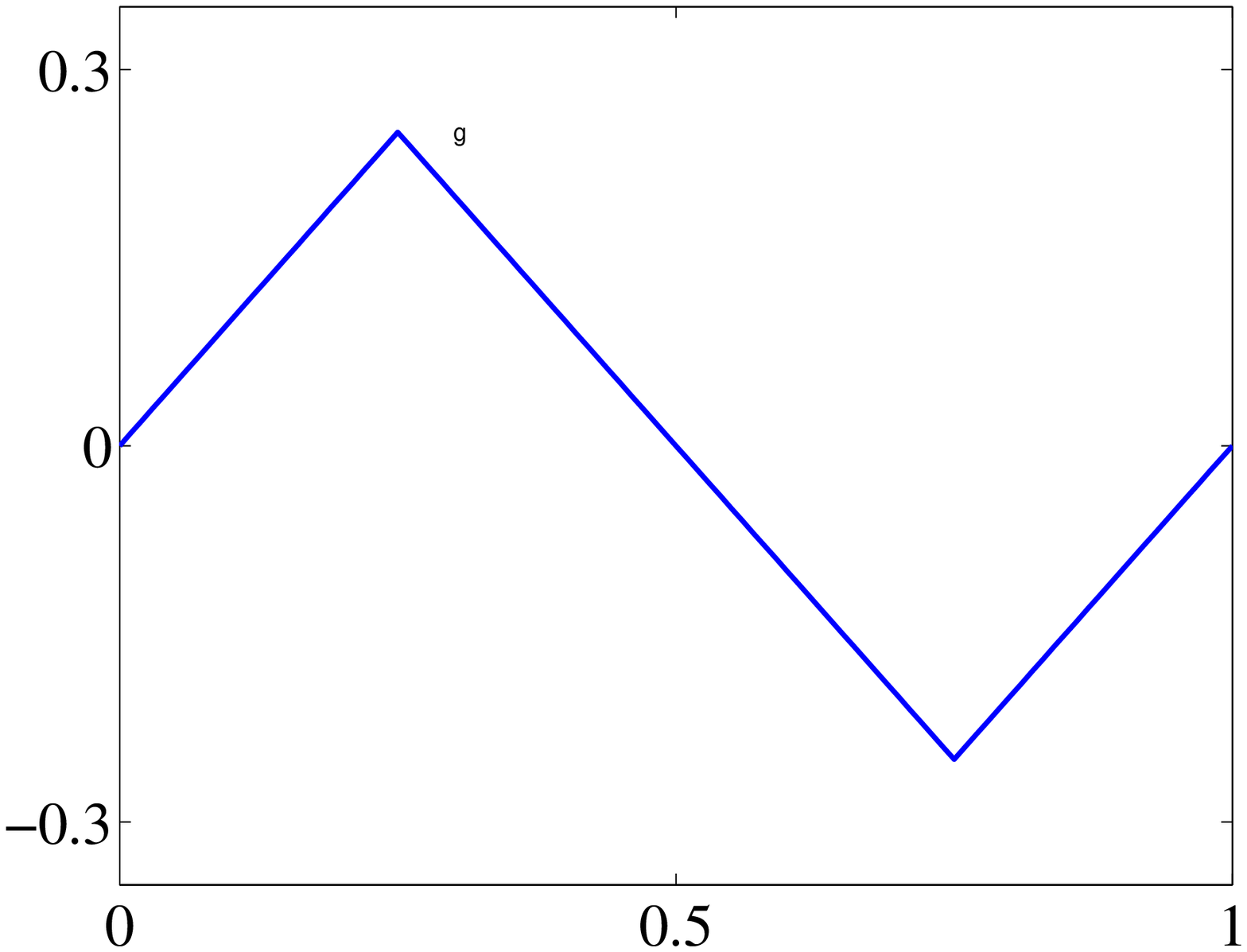} &
      \psfrag{g}{$g_2$}
      \includegraphics[height=.32\columnwidth]{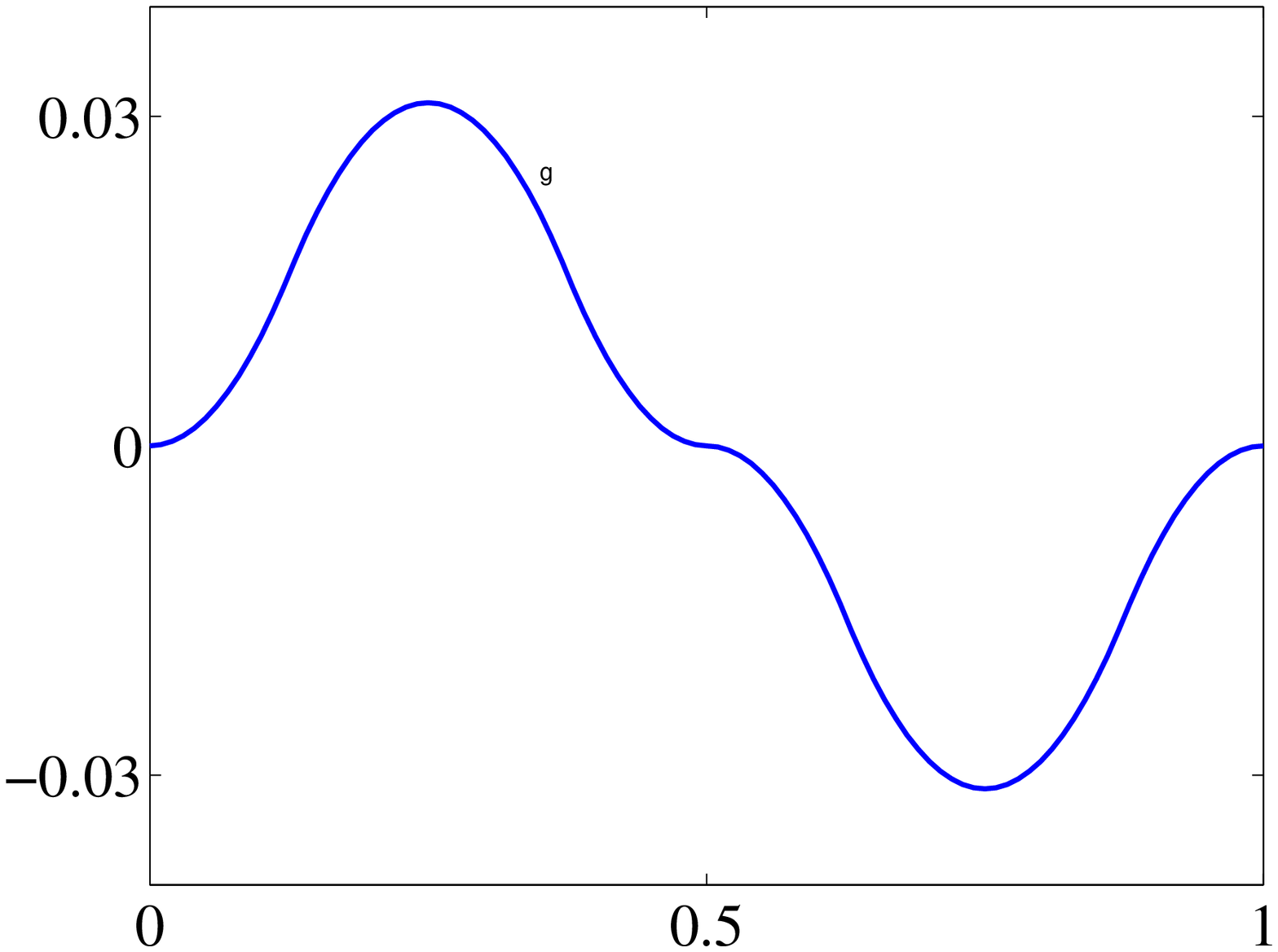} \\
        (a) & (b)
    \end{tabular}
    \caption{\label{fig:bump-func} Panel (a): illustration of
      $1$-Lipschitz continuous bump function $g_1$ used to pack
      $\densclass$ when $\numderiv = 1$.  Panel (b): bump function
      $g_2$ with $|g_2''(\statsample)| \le 1$ used to pack
      $\densclass$ when $\numderiv = 2$.
    }
  \end{center}
\end{figure}

\paragraph{Constructing a good packing:}
We begin by describing the collection of functions we use to prove our lower
bound. Our construction and identification of density functions by vectors is
essentially standard~\cite{Tsybakov09}, but we specify some necessary
conditions that we use later.  First, let $g_\numderiv$ be a function defined
on $[0, 1]$ satisfying the following properties:
\begin{enumerate}[(a)]

\item The function $g_\numderiv$ is $\numderiv$-times differentiable,
  and
\begin{equation*}
  0 = g_\numderiv^{(i)}(0) = g_\numderiv^{(i)}(1/2)
  = g_\numderiv^{(i)}(1)
  ~~ \mbox{for~all~} i < \numderiv.
\end{equation*}

\item The function $g_\numderiv$ is centered with $\int_0^1
  g_\numderiv(x) dx = 0$, and there exist constants $c, c_{1/2} > 0$
  such that
\begin{equation*}
  \int_0^{1/2} g_\numderiv(x) dx = -\int_{1/2}^1 g_\numderiv(x) dx =
  c_{1/2} ~~~ \mbox{and} ~~~ \int_0^1
  \left(g^{(i)}_\numderiv(x)\right)^2 dx \ge c ~~ \mbox{for~all~} i <
  \numderiv.
\end{equation*}
\item The function $g_\numderiv$ is non-negative on $[0, 1/2]$ and
  non-positive on $[1/2, 1]$, and Lebesgue measure is absolutely
  continuous with respect to the measures $G_j, j = 1, 2$ given by
  \begin{equation}
    \label{eqn:absolute-continuity-bumps}
    G_1(A) = \int_{A \cap [0, 1/2]} g_\numderiv(x) dx ~~~ \mbox{and} ~~~
    G_2(A) = -\int_{A \cap [1/2, 1]} g_\numderiv(x) dx.
  \end{equation}
\item Lastly, for almost every $x \in [0, 1]$, we have
  $|g^{(\numderiv)}_\numderiv(x)| \le 1$ and $|g_\numderiv(x)| \le 1$.
\end{enumerate}
\noindent The functions $g_\numderiv$ are smooth ``bumps'' that we use
as pieces in our general construction; see Figure~\ref{fig:bump-func}
for an illustration of such functions in the cases $\numderiv = 1$ and
$\numderiv = 2$

Fix a positive integer $\numbin$ (to be specified momentarily).  Our proof
makes use of the following result from our previous paper~\cite[Lemma
  7]{DuchiJoWa13_parametric}:
\begin{lemma}[Re-stated from the paper~\cite{DuchiJoWa13_parametric}]
  \label{LemPrevious}
  There exists a packing $\packset$ of size at least $\exp(c_0 \numbin)$
  of the hypercube $\{-1, 1\}^\numbin$ such that
  \begin{equation*}
    \lone{\packval - \altpackval} \ge c_1 \numbin \quad \mbox{for all
      $\packval \neq \altpackval$, and}
    \quad
    \frac{1}{|\packset|}
    \sum_{\packval \in \packset} \packval \packval^\top \preceq c_2
    I_{\numbin \times \numbin},
  \end{equation*}
  where $(c_0, c_1, c_2)$ are universal positive constants.
\end{lemma}

We now make use of this packing of the hypercube in order to construct
a packing of our density class. For each $j \in \{1, \ldots,
\numbin\}$, define the function
\begin{align*}
  g_{\numderiv, j}(x) & \defeq \frac{1}{\numbin^\numderiv} \, g_\beta
  \left(\numbin \Big( x - \frac{j - 1}{\numbin}\Big)\right) \indic{ x
    \in [\frac{j-1}{\numbin}, \frac{j}{\numbin}]}.
\end{align*}
Based on this definition, we define the family of
densities
\begin{equation}
  \label{eqn:density-packer}
  \bigg\{ f_\packval \defeq 1 + \sum_{j=1}^\numbin \packval_j
  g_{\numderiv, j} ~~~ \mbox{for}~ \packval \in \packset \bigg\}
  \; \subseteq \densclass.
\end{equation}
It is a standard fact~\cite{Yu97,Tsybakov09} that for any $\packval
\in \packset$, the function $f_\packval$ is $\numderiv$-times
differentiable, satisfies $|f^{(\numderiv)}(x)| \le 1$ for all $x$,
and $\ltwo{f_\packval - f_{\altpackval}}^2 \ge c \numbin^{-2
  \numderiv}$.  Consequently, the class~\eqref{eqn:density-packer} is
a $(c \numbin^{-\numderiv})$-packing of $\densclass$ with cardinality
at least $\exp(c_0 \numbin)$.

\paragraph{Controlling the operator norm of the packing:}
Having constructed a suitable packing of the space $\densclass$, we
now turn to bounding the mutual information associated with a certain
multi-way hypothesis testing problem.  Suppose that an index $\packrv$
is drawn uniformly at random from $\packset$, and conditional on
$\packrv = \packval$, the data points $\statrv_i$ are drawn
i.i.d.\ according to the density $f_\packval$.  The data
$\{\statrv_1, \ldots, \statrv_\numobs \}$ are then passed through an
$\diffp$-locally private distribution $\channelprob$, yielding the
perturbed quantities $\{\channelrv_1, \ldots, \channelrv_\numobs \}$.
The following lemma bounds the mutual information between the
random index $\packrv$ and the outputs $\channelrv_i$.

\begin{lemma}
  \label{lemma:density-information}
  There exists a universal constant $c$ such that for any
  $\diffp$-locally private~\eqref{eqn:local-privacy} conditional
  distribution $\channelprob$, the mutual information is upper bounded
  as
  \begin{equation*}
    \information(\channelrv_1, \ldots, \channelrv_n; \packrv) \le n
    \frac{c \diffp^2}{\numbin^{2 \numderiv + 1}}.
  \end{equation*}
\end{lemma}
\noindent The proof of this claim is fairly involved, so we defer
it to Appendix~\ref{appendix:density-information}. We remark,
however, that standard mutual information bounds~\cite{Yu97,Tsybakov09}
show $\information(\channelrv_1, \ldots, \channelrv_n;
\packrv) \lesssim n / \numbin^{2 \numderiv}$; our bound is thus
essentially a factor of the ``dimension'' $\numbin$ tighter.


\paragraph{Applying testing inequalities:}

The remainder of the proof is an application of Fano's
inequality. In particular, we apply
Lemma~\ref{lemma:minimax-risk-bound} with our $\numbin^{-\numderiv}$
packing of $\densclass$ in $\ltwo{\cdot}$ of size $\exp(c_0 \numbin)$,
and we find that for any $\diffp$-locally private channel $\channelprob$,
there are universal constants $c_0, c_1, c_2$ such that
\begin{equation*}
\minimax_n \left(\densclass, \ltwo{\cdot}^2, \channelprob\right) \ge
\frac{c_0}{\numbin^{2 \numderiv}} \left(1 -
\frac{\information(\channelrv_{1:n}; \packrv) + \log 2}{ c_1 \numbin}
\right) \ge \frac{c_0}{\numbin^{2 \numderiv}} \left(1 - \frac{c_2 n
  \diffp^2 \numbin^{-2 \numderiv - 1} + \log 2}{ c_1 \numbin} \right).
\end{equation*}
Choosing $\numbin_{n, \diffp, \numderiv} = \left(2 c_2 n \diffp^2
\right )^{\frac{1}{2 \numderiv + 2}}$ ensures that the quantity inside
the parentheses is a strictly positive constant. As a consequence,
there are universal constants $c, c' > 0$ such that
\begin{equation*}
  \minimax_n\left(\densclass, \ltwo{\cdot}^2, \diffp\right) \ge
  \frac{c}{\numbin_{n, \diffp, \numderiv}^{2 \numderiv}} = c' \left(n
  \diffp^2\right)^{-\frac{2 \numderiv}{2 \numderiv + 2}}
\end{equation*}
as claimed.


\subsection{Proof of Proposition~\ref{proposition:histogram-estimator}}
\label{sec:proof-histogram}

Note that the operator $\Pi_\numbin$ performs a Euclidean projection
of the vector $(\numbin/n) \sum_{i=1}^n \channelrv_i$ onto the scaled
probability simplex, thus projecting $\what{f}$ onto the set of
probability densities.  Given the non-expansivity of Euclidean
projection, this operation can only decrease the error
$\ltwos{\what{f} - f}^2$. Consequently, it suffices to bound the error
of the unprojected estimator; to reduce notational overhead
we retain our previous notation $\what{\optvar}$ for the unprojected
version.  Using this notation, we have
\begin{align*}
  \E \left [\ltwobig{\what{f} - f}^2\right] & \leq \sum_{j=1}^\numbin
  \E_f\left [\int_{\frac{j-1}{\numbin}}^{\frac{j}{\numbin}}
    (f(\statsample) - \what{\theta}_j)^2 d \statsample \right].
\end{align*}
By expanding this expression and noting that the independent noise
variables $W_{ij} \sim \laplace(\diffp/2)$ have zero mean, we obtain
\begin{align}
  \E \left[\ltwobig{\what{f} - f}^2\right] & \le \sum_{j=1}^\numbin \E_f
  \left[\int_{\frac{j-1}{\numbin}}^{\frac{j}{\numbin}} \bigg
    (f(\statsample) - \frac{\numbin}{n} \sum_{i=1}^n
    [\histelement_\numbin(\statrv_i)]_j\bigg)^2 d\statsample\right] +
  \sum_{j=1}^\numbin \int_{\frac{j-1}{\numbin}}^{\frac{j}{\numbin}}
  \E\bigg[\bigg(\frac{\numbin}{n} \sum_{i = 1}^n W_{ij}\bigg)^2 \bigg]
  \nonumber \\ 
  \label{eqn:private-histogram-risk}
  & = \sum_{j=1}^\numbin \int_{\frac{j-1}{\numbin}}^{\frac{j}{\numbin}}
  \E_f\left[ \bigg(f(\statsample) - \frac{\numbin}{n} \sum_{i=1}^n
    [\histelement_\numbin(\statrv_i)]_j\bigg)^2\right] d\statsample +
  \numbin \, \frac{1}{\numbin} \, \frac{4\numbin^2}{n \diffp^2}.
\end{align}

We bound the error term inside the
expectation~\eqref{eqn:private-histogram-risk}.  Defining $p_j \defeq
\P_f(\statrv \in \statdomain_j) = \int_{\statdomain_j} f(\statsample)
d\statsample$, we have
\begin{equation*}
\numbin\E_f\left[[\histelement_\numbin(\statrv)]_j\right] = \numbin
p_j = \numbin \int_{\statdomain_j} f(\statsample) d\statsample \in
\left[f\left(\statsample\right) - \frac{1}{\numbin},
  f\left(\statsample\right) + \frac{1}{\numbin}\right] ~~
\mbox{for~any~} \statsample \in \statdomain_j,
\end{equation*}
by the Lipschitz continuity of $f$.
Thus, expanding the bias and variance of the integrated expectation above,
we find that
\begin{align*}
  \E_f\left[
    \bigg(f(\statsample) - \frac{\numbin}{n} \sum_{i=1}^n
         [\histelement_\numbin(\statrv_i)]_j\bigg)^2\right]
  & \le \frac{1}{\numbin^2} + \var\left(\frac{\numbin}{n} \sum_{i=1}^n
  \left[\histelement_\numbin(\statrv_i)\right]_j\right) \\
  & = \frac{1}{\numbin^2} + \frac{\numbin^2}{n}
  \var([\histelement_\numbin(\statrv)]_j)
  = \frac{1}{\numbin^2} + \frac{\numbin^2}{n} p_j(1 - p_j).
\end{align*}
Recalling the inequality~\eqref{eqn:private-histogram-risk}, we obtain
\begin{equation*}
  \E_f\left[\ltwobig{\what{f} - f}^2\right]
  \le \sum_{j=1}^\numbin \int_{\frac{j-1}{\numbin}}^\frac{j}{\numbin}
  \left(\frac{1}{\numbin^2} + \frac{\numbin^2}{n} p_j(1 - p_j)\right)
  d \statsample
  + \frac{4\numbin^2}{n \diffp^2}
  = \frac{1}{\numbin^2} + \frac{4\numbin^2}{n\diffp^2}
  + \frac{\numbin}{n} \sum_{j=1}^\numbin p_j(1 - p_j).
\end{equation*}
Since $\sum_{j=1}^\numbin p_j = 1$, we find that
\begin{equation*}
  \E_f\left[\ltwobig{\what{f} - f}^2\right] \le \frac{1}{\numbin^2} +
  \frac{4 \numbin^2}{n \diffp^2} + \frac{\numbin}{n},
\end{equation*}
and choosing $\numbin = (n \diffp^2)^{\frac{1}{4}}$ yields the claim.


\subsection{Proof of Proposition~\ref{proposition:density-upper-bound}}
\label{sec:proof-density-upper-bound}

We begin by fixing $\numbin \in \N$; we will optimize the choice of $\numbin$
shortly.  Recall that, since $f \in \densclass[\lipconst]$, we have $f
= \sum_{j=1}^\infty \optvar_j \basisfunc_j$ for $\optvar_j = \int f
\basisfunc_j$. Thus we may define $\overline{\channelrv}_j = \frac{1}{n}
\sum_{i=1}^n \channelrv_{i,j}$ for each $j \in \{1, \ldots, \numbin\}$, and we have
\begin{equation*}
  \ltwos{\what{f} - f}^2
  = \sum_{j=1}^\numbin (\optvar_j - \overline{\channelrv}_j)^2
  + \sum_{j = \numbin + 1}^\infty \optvar_j^2.
\end{equation*}
Since $f \in \densclass[\lipconst]$, we are guaranteed that
$\sum_{j=1}^\infty j^{2 \numderiv} \optvar_j^2 \leq \lipconst^2$, and
hence
\begin{equation*}
\sum_{j > \numbin} \optvar_j^2 = \sum_{j > \numbin} j^{2 \numderiv}
\frac{\optvar_j^2}{j^{2 \numderiv}} \le \frac{1}{\numbin^{2
    \numderiv}} \sum_{j > \numbin} j^{2 \numderiv} \optvar_j^2 \le
\frac{1}{\numbin^{2 \numderiv}} \lipconst^2.
\end{equation*}
For the indices $j \le \numbin$, we note that by assumption,
$\E[\channelrv_{i,j}] = \int \basisfunc_j f = \optvar_j$,
and since $|\channelrv_{i,j}| \le \sbound$, we have 
\begin{equation*}
  \E\left[(\optvar_j - \overline{Z}_j)^2\right] = \frac{1}{n}
  \var(Z_{1,j}) \le \frac{\sbound^2}{n} =
  \frac{\orthbound^2}{c_\numbin} \, \frac{\numbin}{n} \,
  \left(\frac{e^\diffp + 1}{e^\diffp - 1}\right)^2,
\end{equation*}
where $c_\numbin = \Omega(1)$ is the constant in
expression~\eqref{eqn:size-infinity-channel}.  Putting together the pieces,
the mean-squared $L^2$-error is upper bounded as
\begin{equation*}
\E_f\left[\ltwos{\what{f} - f}^2\right] \le c \left(\frac{\numbin^2}{n
  \diffp^2} + \frac{1}{\numbin^{2 \numderiv}}\right),
\end{equation*}
where $c$ is a constant depending on $\orthbound$, $c_\numbin$, and
$\lipconst$.  Choose $\numbin = (n \diffp^2)^{1 / (2 \numderiv +
  2)}$ to complete the proof.

\subsection{Proof of Observation~\ref{observation:laplace-density-bad}}
\label{sec:proof-laplace-density-bad}

We begin by noting that for $f = \sum_j \optvar_j \basisfunc_j$,
by definition of $\what{f} = \sum_j \what{\optvar}_j \basisfunc_j$ we have
\begin{equation*}
  \E\left[\ltwos{f - \what{f}}^2\right]
  = \sum_{j=1}^\numbin \E\left[(\optvar_j - \what{\optvar}_j)^2\right]
  + \sum_{j \ge \numbin + 1} \optvar_j^2
  = \sum_{j=1}^\numbin \frac{\orthbound^2 \numbin^2}{n \diffp^2}
  + \sum_{j \ge \numbin + 1} \optvar_j^2
  = \frac{\orthbound^2 \numbin^3}{n \diffp^2}
  + \sum_{j \ge \numbin + 1} \optvar_j^2.
\end{equation*}
Without loss of generality, let us assume $\basisfunc_1 = 1$ is the constant
function. Then $\int \basisfunc_j = 0$ for all $j > 1$, and by defining the
true function $f = \basisfunc_1 + (\numbin + 1)^{-\numderiv}
\basisfunc_{\numbin + 1}$, we have $f \in \densclass[2]$ and
$\int f = 1$, and moreover,
\begin{equation*}
  \E\left[\ltwos{f - \what{f}}^2\right]
  \ge \frac{\orthbound^2 \numbin^3}{n \diffp^2}
  + \frac{1}{(\numbin + 1)^{-2 \numderiv}}
  \ge
  C_{\numderiv, \orthbound}
  (n \diffp^2)^{-\frac{2 \numderiv}{2 \numderiv + 3}},
\end{equation*}
where $C_{\numderiv, \orthbound}$ is a constant depending on
$\numderiv$ and $\orthbound$.  This final lower bound comes by
minimizing over all $\numbin$.  (If $(\numbin + 1)^{-\numderiv}
\orthbound > 1$, we can rescale $\basisfunc_{\numbin + 1}$ by
$\orthbound$ to achieve the same result and guarantee that $f \ge 0$.)


\section{Discussion}

We have
linked minimax analysis from statistical decision theory with
differential privacy, bringing some of their respective foundational
principles into close contact.  Our main technique, in the form of the
divergence bounds in Proposition~\ref{proposition:information-bounds},
shows that applying differentially private sampling schemes
essentially acts as a contraction on distributions, and we think that
such results may be more generally applicable. In this paper
particularly, we showed how to apply our divergence bounds to obtain
sharp bounds on the convergence rate for certain nonparametric
problems in addition to standard finite-dimensional settings.  With
our earlier paper~\cite{DuchiJoWa13_parametric}, we have developed a
set of techniques that show that roughly, if one can construct a
family of distributions $\{\statprob_\packval\}$ on the sample space
$\statdomain$ that is not well ``correlated'' with any member of $f
\in L^\infty(\statdomain)$ for which $f(\statsample) \in \{-1, 1\}$,
then providing privacy is costly---the contraction
Proposition~\ref{proposition:information-bounds} provides is strong.

By providing (to our knowledge, the first) sharp convergence rates for many
standard statistical inference procedures under local differential privacy, we
have developed and explored some tools that may be used to better understand
privacy-preserving statistical inference and estimation procedures. We have
identified a fundamental continuum along which privacy may be traded for
utility in the form of accurate statistical estimates, providing a way to
adjust statistical procedures to meet the privacy or utility needs of the
statistician and the population being sampled. Formally identifying this
tradeoff in other statistical problems should allow us to better understand
the costs and benefits of privacy; we believe we
have laid some of the groundwork for doing so.

\subsection*{Acknowledgments}

We thank Guy Rothblum for very helpful discussions.  JCD was supported
by a Facebook Graduate Fellowship and an NDSEG fellowship.  Our work
was supported in part by the U.S.\ Army Research Laboratory,
U.S.\ Army Research Office under grant number W911NF-11-1-0391, and
Office of Naval Research MURI grant N00014-11-1-0688.


\appendix

\section{Proof of Lemma~\ref{lemma:lots-of-packings}}
\label{appendix:proof-lemma-lots-of-packings}

In the regime $\nnz \in \openleft{d/2}{d}$, the statement of lemma follows
Lemma~7 in~\citet{DuchiJoWa13_parametric};
consequently, we prove the claim for $\nnz \leq d/2$.  If
$\nnz \in \{1, 2, 3, 4 \}$, taking $\packset = \{\packval \in H_d :
\lone{\packval} = \nnz\}$ implies that $\lone{\packval -
  \altpackval} \ge \nnz/4$ for $\packval \neq \altpackval$,
$\card(\packset) = \binom{d}{\nnz} \ge (d/\nnz)^{c \nnz}$ for some
constant $c > 0$, and
\begin{equation*}
  \cov(\packrv) = \left(\frac{\nnz}{d} - \frac{\nnz^2}{d^2}\right)
  I_{d \times d},
\end{equation*}
from which the claim follows.

Accordingly, we focus on $\nnz \in \openleft{4}{d/2}$. To further simplify the
analysis, we claim it suffices to establish the claim in the case that
$\nnz/4$ is integral (i.e., $\nnz \in 4 \N$).  Indeed, assume that the result
holds for all such integers. Given some $\nnz \notin 4 \N$, we may consider a
packing $\packset'$ of the binary hypercube $H_{d'}$ with $d' = d - (\nnz - 4
\floor{\nnz/4})$ and $\lone{\packval} = \nnz' = 4\floor{\nnz/4}$ for
$\packval \in \packset'$.
By assumption, there is a packing $\packset'$ of $H_{d'}$
satisfying the lemma. Now to each vector $\packval \in \packset'$, we
concatenate the $(\nnz - 4 \floor{\nnz/4})$-vector $\onevec$, which gives
$[\packval^\top ~ \onevec^\top]^\top \in \{0, 1\}^d$ and $\lone{[\packval^\top
    ~ \onevec^\top]} = \nnz$. This concatenation does not increase
$\cov(\packrv)$---the last $\nnz - 4\floor{\nnz/4}$ coordinates have
covariance zero---and the rest of the terms in items
(\ref{item:nnz})--(\ref{item:sparse-covariance}) incur only constant factor
changes. \\

It remains to prove the claim for $\nnz \in 4 \N$ over the range $\{5,
\ldots, \floor{d/2} \}$.  To ease notation, we let $\locnnz = \nnz/4$
belong to the interval $[2, d/8]$.  Our proof is based on the
probabilistic method~\cite{AlonSp00}: we propose a random construction
of a packing, and show that it satisfies the desired properties with
strictly positive probability.  Our random construction is
straightforward: letting $H_d = \{0,1\}^d$ denote the Boolean
hypercube, we sample $\Mpack$ i.i.d.\ random vectors $U_i$ from the uniform
distribution over the set
\begin{equation}
  \label{eqn:samplable-set}
  S_\locnnz \defeq \left \{\packval \in H_d \, \mid \, \lone{\packval}
  \defeq 4 \locnnz \right \}.
\end{equation}
We claim that for $\Mpack = (d/(6\locnnz))^{3\locnnz/2}$, the
resulting random set $\USET \defeq \{U_1, \ldots, U_\Mpack\}$
satisfies the claimed properties with non-zero probability. We say that
$\USET$ is $\locnnz$-separated if $\|U_i -
U_j\|_1 > \locnnz$ for all $i \neq j$, and we use $\cov(\USET)$ to
denote the covariance of a random vector $V$ drawn uniformly at random
from $\USET$.  Our proof is based on the following two tail bounds,
which we prove shortly: for a universal constant $c < \infty$,
\begin{subequations}
  \begin{align}
    \label{EqnGoodSep}
    \P \left[ \mbox{$\USET$ is not $\locnnz$-separated}\right]
    & \leq \binom{\Mpack}{2}
    \left(6 \, \frac{\locnnz}{d}\right)^{3\locnnz}, \qquad
    \mbox{and} \\
    \label{EqnGoodCov}
    \P \left[\lambda_{\max} \big(\cov(\USET) \big) \geq t \right] & 
    \le d
    \exp \left(-\frac{\Mpack t^2}{3c \max\{\locnnz, \locnnz^3/d\} + c t
      \locnnz}\right) \qquad \mbox{for all $t > 0$.}
  \end{align}
\end{subequations}

For the moment, let us assume the validity of these bounds and use
them to complete the proof.  By the union bound, we have
\begin{align*}
  \P\left(\mbox{$\USET$ is not $\locnnz$-separated} \mbox{ or}~
  \cov(\USET) \not \preceq t I\right) & \leq \binom{\Mpack}{2} \left(6
  \, \frac{\locnnz}{d}\right)^{3\locnnz} + d \exp\left(-\frac{\Mpack
    t^2}{3c \max\{\locnnz, \locnnz^3/d\} + c t \locnnz} \right).
\end{align*}
By choosing $t = C \locnnz / d$ and recalling that $\Mpack =
(d/(6\locnnz))^{3\locnnz/2}$, we obtain the bound
\begin{equation*}
  \half + d \exp\left(-C^2 \frac{\locnnz^2 (d / (6\locnnz))^{3\locnnz/2}}{
    3c \max\{d^2 \locnnz, d \locnnz^3\} +  C c d \locnnz^2}\right).
\end{equation*}
If $\locnnz \ge \locnnz^3/d$, the second term can be easily seen to be less
than $\half$ for suitably large constants $C$, so assume that $\locnnz \le
\locnnz^3 / d$.  Then we have, where $c$ is a constant whose value may change
from inequality to inequality,
\begin{equation*}
  \frac{\locnnz^2 (d / (6\locnnz))^{3\locnnz/2}}{
    3c \max\{d^2 \locnnz, d \locnnz^3\} + C c d \locnnz^2}
  =
  \frac{\locnnz^2 d^{3\locnnz/2}}{(6\locnnz)^{3\locnnz/2}(3c d \locnnz^3
    + C c d \locnnz^2)}
  \ge c \frac{d^{3\locnnz/2}}{(6\locnnz)^{3\locnnz/2} d \locnnz}
  \ge \frac{c}{d\locnnz} \left(\frac{d}{6\locnnz}\right)^{3\locnnz / 2}.
\end{equation*}
For suitably large $d$ and any $\locnnz \ge 2$, the final term is greater
than $c' \log d$ for some constant $c' > 0$, which implies that
with appropriate choice of the constant $C$ earlier, we have the
bound
\begin{align*}
\P\left(\mbox{$\USET$ is not $\locnnz$-separated} \mbox{ or}~
\cov(\USET) \not \preceq t I\right) & < 1.
\end{align*}
Consequently, recalling that $\nnz = 4\locnnz$ by definition, a
packing as described in the statement of the lemma must exist. \\

It remains to prove the tail bounds~\eqref{EqnGoodSep}
and~\eqref{EqnGoodCov}.  Beginning with the former bound, define the
set
\begin{equation*}
  N(\packval, \locnnz) \defeq \left\{ \altpackval \in H_d \, \mid \,
  \lone{\packval - \altpackval} \le \locnnz \right\}.
\end{equation*}
Recalling the definition~\eqref{eqn:samplable-set} $S_\locnnz$,
let $U_i$ and $U_j$ be sampled independently and uniformly at random from
$S_\locnnz$. Then
\begin{align*}
  \P \left( \|U_i - U_j\|_1 \leq \locnnz \right) & \leq
  \frac{\card(N(\packval, \locnnz))}{\card(S_\locnnz)}.
\end{align*}
Note that $N(\packval, \locnnz)$ can be constructed by choosing an
arbitrary subset $J \subset \{1, \ldots, d\}$ of size $\locnnz$, and
then setting $\altpackval_j = \packval_j$ for $j \not \in J$ and
$\packval_j$ arbitrarily otherwise; consequently, its cardinality is
upper bounded as $\card(N(\packval, \locnnz)) \leq \binom{d}{\locnnz}
2^\locnnz$.  Since $\card(S_\locnnz) = \binom{d}{4\locnnz}$, we find
that
\begin{equation*}
  \frac{\card(N(\packval, \locnnz))}{\card(S_\locnnz)} =
  \frac{\binom{d}{\locnnz} 2^\locnnz}{\binom{d}{4\locnnz}} =
  \frac{2^\locnnz (d - 4\locnnz)! (4\locnnz)!}{(d - \locnnz)!
    \locnnz!}  = 2^\locnnz \prod_{j=1}^{3\locnnz} \frac{\locnnz + j}{d
    - 4\locnnz + j} \le 2^\locnnz \left(\frac{4\locnnz}{d -
    \locnnz}\right)^{3\locnnz},
\end{equation*}
where the final inequality follows because the function $x \mapsto
h(x) = \frac{\locnnz + x}{d - 4\locnnz + x}$ is increasing for $x >
0$.  Since $\locnnz \le d/8$ by assumption, we arrive at the upper
bound
\begin{equation*}
  \frac{\card(N(\packval, \locnnz))}{\card(S_\locnnz)} \le 2^\locnnz
  \left(\frac{4\locnnz}{d - \locnnz}\right)^{3\locnnz} = \left(\frac{4
    \cdot 2^{1/3} \locnnz}{d - \locnnz}\right)^{3\locnnz} \le \left(6 \,
  \frac{\locnnz}{d}\right)^{3\locnnz}.
\end{equation*}
Since we have to compare $\binom{\Mpack}{2}$ such pairs over the set
$\USET$, the claim~\eqref{EqnGoodSep} follows from the union bound.


We now turn to establishing the claim~\eqref{EqnGoodCov}, for which we
make use of matrix Bernstein inequalities.  Letting $U$ be drawn
uniformly at random from $S_\locnnz$, we have
\begin{equation*}
  \E[U U^\top] = \beta_{\locnnz, d} \onevec \onevec^\top +
  \left(\frac{4 \locnnz}{d} - \beta_{\locnnz, d}\right) I_{d \times d}.
\end{equation*}
where $\beta_{\locnnz,d} \defeq \binom{4 \locnnz}{2} \binom{d}{2}^{-1}$.
Consequently, the $d \times d$ random matrix
\begin{align*}
A & \defeq UU^\top - \beta_{\locnnz,d} \onevec \onevec^\top -
\left(\frac{4 \locnnz}{d} - \beta_{\locnnz,d}\right) I_{d \times d}.
\end{align*}
is centered ($\E[A] = 0$), and by definition of our construction,
$\cov(\USET) = \frac{1}{\Mpack} \sum_{i=1}^\Mpack A_i$, where the
random matrices $\{A_i\}_{i=1}^\Mpack$ are drawn i.i.d.

In order to apply a matrix Bernstein inequality, it remains to bound
the operator norm (maximum singular value) of $A$ and its variance.
The operator norm of $A$ is upper bounded as
\begin{equation*}
  \matrixnorm{A} \le \matrixnorm{UU^\top - (4\locnnz/d - \beta_{\locnnz,d})
    I} + \beta_{\locnnz,d} \matrixnorm{\onevec \onevec^\top}
  = 4\locnnz - \frac{4 \locnnz}{d}
  + \beta_{\locnnz, d} + d \beta_{\locnnz,d}
  \le 5 \locnnz.
\end{equation*}
Moreover, we claim that there is a universal positive constant $c$
such that 
\begin{align}
  \label{EqnIntermediate}
  \matrixnorm{\E[A^2]} & \leq c \max\{\locnnz, \locnnz^3 / d\}.
\end{align}
To establish this claim, we begin by computing
\begin{align*}
  \E[A^2] & = \E[UU^\top UU^\top] - \left(\left(\frac{4 \locnnz}{d} -
  \beta_{\locnnz,d}\right) I_{d \times d} + \beta_{\locnnz,d} \onevec
  \onevec^\top\right)^2 \\
  & = 4 \locnnz \left(\left(\frac{4 \locnnz}{d} -
  \beta_{\locnnz,d}\right) I_{d \times d} + \beta_{\locnnz,d} \onevec
  \onevec^\top\right) - \left(\left(\frac{4 \locnnz}{d} -
  \beta_{\locnnz,d}\right) I_{d \times d} + \beta_{\locnnz,d} \onevec
  \onevec^\top\right)^2.
\end{align*}
Consequently, if we define the constants,
\begin{equation*}
  a_{\locnnz, d} \defeq \left(4 \locnnz - \frac{4 \locnnz}{d} +
  \beta_{\locnnz,d}\right) ~~~ \mbox{and} ~~~ b_{\locnnz, d} \defeq
  \left(4 \locnnz \beta_{\locnnz, d} - \frac{8 \locnnz
    \beta_{\locnnz,d}}{d} + 2 \beta_{\locnnz,d}^2 - d \beta_{\locnnz,
    d}^2\right),
\end{equation*}
then $\E[A^2] = a_{\locnnz,d}I_{d \times d} + b_{\locnnz, d} \onevec
\onevec^\top$.  It is easy to see that $|a_{\locnnz,d}| \le 4 \locnnz$
and that $|b_{\locnnz,d}| \le c' \frac{\locnnz^3}{d^2}$ for some
universal constant $c'$, from which the intermediate
claim~\eqref{EqnIntermediate} follows.  With these pieces in place,
the claimed tail bound~\eqref{EqnGoodCov} follows a matrix Bernstein
inequality (e.g.,~\cite[Corollary 5.2]{MackeyJoChFaTr12}), applied to
the quantity $\cov(\USET) = \frac{1}{\Mpack} \sum_{i=1}^\Mpack A_i$.


\section{Proof of Lemma~\ref{lemma:density-information}}
\label{appendix:density-information}

This result relies on inequality~\eqref{eqn:super-fano} from
Proposition~\ref{proposition:information-bounds}, along with a careful
argument to understand the extreme points of $\optdens \in
L^\infty([0, 1])$ that we use when applying the proposition. First,
we take a packing $\packset$ as guaranteed by Lemma~\ref{LemPrevious}
and consider densities $f_\packval$ for $\packval \in \packset$.
Overall,
our first step is to show for the purposes of applying
inequality~\eqref{eqn:super-fano}, it is no loss of generality to
identify $\optdens \in L^\infty([0, 1])$ with vectors $\optdens \in
\R^{2 \numbin}$, where $\optdens$ is constant on intervals of the form
$[i/2\numbin, (i + 1)/2\numbin]$. With this identification complete,
we can then use the packing set $\packset$ from
Lemma~\ref{LemPrevious} to provide a bound on the correlation of
any $\optdens \in L^\infty([0, 1])$ with the densities $f_\packval$,
which completes the proof.

With this outline in mind, let the sets $\MySet_i$, $i \in \{1, 2,
\ldots, 2\numbin\}$, be defined as $\MySet_i =
\openright{(i-1)/2\numbin}{i/2\numbin}$ except that $\MySet_{2\numbin}
= [(2\numbin - 1)/2\numbin, 1]$, so the collection
$\{\MySet_i\}_{i=1}^{2 \numbin}$ forms a partition of the unit
interval $[0, 1]$. By construction of the densities $f_\packval$, the
sign of $f_\packval - 1$ remains constant on each
$\MySet_i$. Recalling the linear functionals $\varphi_\packval$ in
Proposition~\ref{proposition:information-bounds}, we have
$\varphi_\packval : L^\infty([0, 1]) \rightarrow \R$ defined via
\begin{equation*}
  \varphi_\packval(\optdens) = \sum_{i=1}^{2\numbin} \int_{\MySet_i}
  \optdens(x)( f_\packval(x) - \overline{f}(x)) dx, =
  \sum_{i=1}^{2\numbin} \int_{\MySet_i} \optdens(x)(f_\packval(x) - 1
  - (\overline{f}(x) - 1)) dx,
\end{equation*}
where $\overline{f} = (1 / |\packset|) \sum_{\packval \in \packset}
f_\packval$. Expanding the square, we find that
since $\overline{f}$ is the average, we have
\begin{equation*}
  \frac{1}{|\packset|}
  \sum_{\packval \in \packset}
  \varphi_\packval(\optdens)^2
  \le \frac{1}{|\packset|}
  \sum_{\packval \in \packset}
  \bigg(\sum_{i=1}^{2\numbin} \int_{\MySet_i} \optdens(x)(f_\packval(x) - 1) dx
  \bigg)^2.
\end{equation*}

Since the set $\linfset$ from
Proposition~\ref{proposition:information-bounds} is compact, convex,
and Hausdorff, the Krein-Milman theorem~\cite[Proposition
  1.2]{Phelps01} guarantees that it is equal to the convex hull of its
extreme points; moreover, since the functionals $\optdens \mapsto
\varphi^2_\packval(\optdens)$ are convex, the supremum in
Proposition~\ref{proposition:information-bounds} must be attained at
the extreme points of $\linfset$.  As a consequence, when applying the
information bound
\begin{equation}
  \information(\channelrv_1, \ldots, \channelrv_n; \packrv) \le n
  C_\diffp \frac{1}{|\packset|} \sup_{\optdens \in \linfset}
  \sum_{\packval \in \packset} \varphi_\packval^2(\optdens),
  \label{eqn:information-sup-densities}
\end{equation}
we can restrict our attention to $\optdens \in L^\infty([0, 1])$ for
which $\optdens(x) \in \{e^{-\diffp} - e^\diffp, e^\diffp -
e^{-\diffp}\}/2$.

Now we argue that it is no loss of generality to assume that
$\optdens$, when restricted to $\MySet_i$, is a constant (apart from a
measure zero set). Using $\lebesgue$ to denote Lebesgue measure,
define the shorthand $\kappa = (e^{\diffp} - e^{-\diffp}) / 2$. Fix $i
\in [2\numbin]$, and assume for the sake of contradiction that there
exist sets $B_i, C_i \subset \MySet_i$ such that $\optdens(B_i) =
\{\kappa\}$ and $\optdens(C_i) = \{-\kappa\}$, while $\lebesgue(B_i) >
0$ and $\lebesgue(C_i) > 0$.\footnote{For a function
  $f$ and set $A$, the notation $f(A)$
  denotes the image $f(A) = \{f(x) \mid x \in A\}$.}
We will construct vectors $\optdens_1$
and $\optdens_2 \in \linfset$ and a value $\lambda \in (0, 1)$ such
that
\begin{equation*}
  \int_{\MySet_i} \optdens(x) (f_\packval(x) - 1)
  d x =
  \lambda \int_{\MySet_i} \optdens_1(x) (f_\packval(x) - 1) dx
  + (1 - \lambda) \int_{\MySet_i} \optdens_2(x) (f_\packval(x) - 1) dx
\end{equation*}
simultaneously for all $\packval \in \packset$, while
on $\MySet_i^c = [0, 1] \setminus \MySet_i$, we will have the equivalence
\begin{equation*}
  \left.\optdens_1\right|_{\MySet_i^c}
  \equiv \left.\optdens_2 \right|_{\MySet_i^c}
  \equiv \left.\optdens\right|_{\MySet_i^c}.
\end{equation*}
Indeed, set $\optdens_1(\MySet_i) = \{\kappa\}$ and $\optdens_2(\MySet_i) =
\{-\kappa\}$, otherwise setting $\optdens_1(x) = \optdens_2(x) = \optdens(x)$
for $x \not \in \MySet_i$.  We define
\begin{equation*}
  \lambda \defeq \frac{\int_{B_i} (f_\packval(x) - 1) dx}{
    \int_{\MySet_i} (f_\packval(x) - 1) dx}
  ~~~ \mbox{so} ~~~
  1 - \lambda = \frac{\int_{C_i} (f_\packval(x) - 1) dx}{
    \int_{\MySet_i} (f_\packval(x) - 1) dx}.
\end{equation*}
By the construction of the function $g_\numderiv$, the function
$f_\packval - 1$ does not change signs on $\MySet_i$, and the absolute
continuity conditions on $g_\numderiv$ specified in
equation~\eqref{eqn:absolute-continuity-bumps} guarantee $1 > \lambda
> 0$, since $\lebesgue(B_i) > 0$ and $\lebesgue(C_i) > 0$.  Moreover,
the quantity $\lambda$ is constant for all $\packval$ by the
construction of the $f_\packval$, since $B_i \subset \MySet_i$ and $C_i
\subset \MySet_i$.  We thus find that for any $\packval \in \packset$,
\begin{align*}
  \lefteqn{\int_{\MySet_i} \optdens(x)  (f_\packval(x) - 1) dx
    = \int_{B_i} \optdens_1(x)
    (f_\packval(x) - 1)dx
    + \int_{C_i} \optdens_2(x) (f_\packval(x) - 1) dx} \\ 
  & = \kappa \int_{B_i} (f_\packval(x) - 1)dx
  - \kappa \int_{C_i} (f_\packval(x) - 1) dx
  = \kappa \lambda \int_{\MySet_i} (f_\packval(x) - 1)dx
  - \kappa (1 - \lambda) \int_{\MySet_i} (f_\packval(x) - 1)dx \\
  & = \lambda \int \optdens_1(x)
  (f_\packval(x) - 1)dx + (1 - \lambda) \int \optdens_2(x)
  (f_\packval(x) - 1)dx.
\end{align*}
By linearity and the strong convexity of the function $x \mapsto x^2$,
then, we find that
\begin{align*}
  \lefteqn{\sum_{\packval \in \packset}
    \bigg(\sum_{i=1}^{2\numbin} \int_{\MySet_i}
    \optdens(x) (f_\packval(x) - 1)dx \bigg)^2} \\
  & < \lambda \sum_{\packval
    \in \packset} \bigg(\sum_{i=1}^{2\numbin} \int_{\MySet_i} \optdens_1(x)
  (f_\packval(x) - 1) dx\bigg)^2
  + (1 - \lambda) \sum_{\packval \in
    \packset} \bigg(\sum_{i=1}^{2\numbin} \int_{\MySet_i} \optdens_2(x)
  (f_\packval(x) - 1)dx \bigg)^2.
\end{align*}
Thus one of the densities $\optdens_i$, $i \in \{1, 2\}$ must have a
larger objective value than $\optdens$. This is our desired
contradiction, which shows that (up to
measure zero sets) any $\optdens$ attaining the supremum in the
information bound~\eqref{eqn:information-sup-densities} must be
constant on each of the $\MySet_i$.

Having shown that $\optdens$ is constant on each of the intervals
$\MySet_i$, we conclude that the
supremum~\eqref{eqn:information-sup-densities} can be reduced to a
finite-dimensional problem over the subset
\begin{align*}
  \linfsetfinite{2\numbin} \defeq \left\{u \in \R^{2\numbin}
  \mid \linf{u} \le \frac{e^\diffp
    - e^{-\diffp}}{2} \right\}
\end{align*}
of $\R^{2\numbin}$.
In terms of this subset, we have the upper bound
\begin{equation*}
\frac{|\packset|}{C_\diffp n} \, \information(\channelrv_1, \ldots,
\channelrv_n; \packrv) \le \sup_{\optdens \in \linfset} \sum_{\packval
  \in \packset} \varphi_\packval(\optdens)^2 \le \sup_{\optdens \in
  \linfsetfinite{2 \numbin}}
  \sum_{\packval \in \packset}
  \bigg(\sum_{i=1}^{2\numbin} \optdens_i
  \int_{\MySet_i} (f_\packval(x) - 1) dx\bigg)^2.
\end{equation*}
By construction of the $f_\packval$ and $g_\numderiv$, we have the
equality
\begin{equation*}
\int_{\MySet_i} (f_\packval(x) - 1) dx = (-1)^{i + 1}\packval_i
\int_{0}^{\frac{1}{2\numbin}} g_{\numderiv, 1}(x) dx = (-1)^{i + 1}
\packval_i \int_0^{\frac{1}{2\numbin}} \frac{1}{\numbin^\numderiv}
g(\numbin x) dx = (-1)^{i + 1} \packval_i
\frac{c_{1/2}}{\numbin^{\numderiv + 1}},
\end{equation*}
which implies that
\begin{align}
\frac{|\packset|}{C_\diffp \numobs} \information(\channelrv_1, \ldots,
\channelrv_n; \packrv) & \le \sup_{\optdens \in \linfsetfinite{2
    \numbin}} \sum_{\packval \in \packset}
\left(\frac{c_{1/2}}{\numbin^{\numderiv + 1}} \optdens^\top \left(
\left[\begin{matrix} 1 \\ -1 \end{matrix}\right] \otimes \packval
\right)\right)^2 \nonumber \\ & = \frac{c_{1/2}^2}{\numbin^{2
    \numderiv + 2}} \sup_{\optdens \in \linfsetfinite{2 \numbin}}
\optdens^\top \bigg(\sum_{\packval \in \packset} \left[\begin{matrix}
    1 \\ -1 \end{matrix}\right] \otimes \packval \packval^\top \otimes
\left[\begin{matrix} 1 \\ -1 \end{matrix}\right]^\top \bigg)\optdens,
  \label{eqn:kronecker-density-bound}
\end{align}
where $\otimes$ denotes the Kronecker product.  By our construction of
the packing $\packset$ of $\{-1, 1\}^\numbin$, there exists a constant
$c$ such that $(1/|\packset|) \sum_{\packval \in \packset} \packval
\packval^\top \preceq c I_{\numbin \times \numbin}$.  Moreover,
observe that the mapping
\begin{equation*}
  A \mapsto \left[\begin{matrix} 1 \\ -1 \end{matrix}\right]
  \otimes A \otimes \left[\begin{matrix} 1
      \\ -1 \end{matrix}\right]^\top ~~~ \mbox{satisfies} ~~~
  \left[\begin{matrix} x \\ y \end{matrix}\right]^\top
  \bigg(\left[\begin{matrix} 1 \\ -1 \end{matrix}\right] \otimes A
  \otimes \left[\begin{matrix} 1 \\ -1 \end{matrix}\right]^\top \bigg)
  \left[\begin{matrix} x \\ y \end{matrix}\right] = (x - y)^\top A (x
  - y),
\end{equation*}
whence it is operator monotone ($A \succeq B$ implies $(x - y)^\top A
(x - y) \ge (x - y)^\top B (x - y)$).  Consequently, by linearity of the
Kronecker product $\otimes$ and Lemma~\ref{LemPrevious}, there
is a universal constant $c$ such that
\begin{equation*}
  \frac{1}{|\packset|} \sum_{\packval \in \packset}
  \left[\begin{matrix} 1 \\ -1 \end{matrix}\right] \otimes \packval
  \packval^\top \otimes \left[\begin{matrix} 1
      \\ -1 \end{matrix}\right]^\top = \left[\begin{matrix} 1
      \\ -1 \end{matrix}\right] \otimes \bigg(\frac{1}{|\packset|}
  \sum_{\packval \in \packset} \packval \packval^\top \bigg) \otimes
  \left[\begin{matrix} 1 \\ -1 \end{matrix}\right]^\top \preceq c
  I_{2\numbin \times 2\numbin}.
\end{equation*}
Combining this bound with our
inequality~\eqref{eqn:kronecker-density-bound}, we obtain
\begin{align*}
\information(\channelrv_1, \ldots, \channelrv_\numobs; \packrv) \le
\numobs \frac{c}{\numbin^{2 \numderiv + 2}} \sup_{\optdens \in
  \linfsetfinite{2 \numbin}} \optdens^\top I \optdens = \numobs
\frac{c (e^\diffp - e^{-\diffp})^2 \numbin}{2 \numbin^{2 \numderiv +
    2}} 
\end{align*}
for some universal numerical constant $c$.  Since $\diffp \in [0, 1/4]$,
we have $(e^\diffp - e^{-\diffp})^2 \leq c' \diffp^2$, which completes
the proof.


\bibliographystyle{abbrvnat} \bibliography{bib}

\end{document}